\documentclass[11pt, a4paper]{article}

\usepackage[left=2.5cm,right=2.5cm,top=2.5cm,bottom=2cm]{geometry}
\usepackage[english]{babel}
\usepackage{parskip}
\usepackage[utf8]{inputenc}
\usepackage[colorinlistoftodos]{todonotes}
\usepackage{graphicx}
\usepackage{titlesec}
\usepackage{pdfpages}
\usepackage[version=4]{mhchem}
\usepackage{amsfonts}
\usepackage[font=small]{caption}
\usepackage{amsmath}
\usepackage{mathtools}
\usepackage{nicefrac}
\usepackage[thinc]{esdiff}
\usepackage{derivative}
\usepackage{siunitx}
\usepackage{array}
\usepackage{booktabs}
\usepackage{multirow}
\usepackage[affil-it]{authblk}
\usepackage[colorinlistoftodos]{todonotes}
\usepackage{hyperref}
\usepackage[doi=false,isbn=false,url=false,eprint=false,backend=biber,style=numeric-comp, giveninits, sorting=none,maxnames=50]{biblatex}

\DeclareSIUnit\year{yr}

\usepackage{rotating}					
\usepackage{longtable}					
\usepackage{parcolumns}					
\usepackage{tabularx} 					
\usepackage{makecell}                   
\usepackage{placeins}                   
\usepackage{float}  
\usepackage{textcomp}                   
\usepackage{threeparttable}
\usepackage{caption}
\DeclareCaptionLabelFormat{AppendixTables}{A.#2} 
\usepackage{subcaption}
\usepackage{csquotes}
\usepackage[title]{appendix}

\usepackage{setspace} 
\usepackage{lineno}	 

\titleformat{\section}
{\large\bfseries}
{\thesection.}{0.5em}{}
\titlespacing{\section}{0mm}{1pc}{0.8pc}

\titleformat{\subsection}
{\normalfont\bfseries}
{\thesubsection.}{0.5em}{}
\titlespacing{\subsection}{0cm}{0.8pc}{0.4pc}

\titleformat{\subsubsection}
{\normalfont\bfseries}
{\thesubsubsection.}{0.5em}{}
\titlespacing{\subsubsection}{0cm}{0.8pc}{0.4pc}

\graphicspath{{Pictures/}}
\pretocmd{\section}{\FloatBarrier}{}{} 
\title{Cost-Optimal Power-to-Methanol: Flexible Operation or Intermediate Storage?}

\author{Simone Mucci $^1$, Alexander Mitsos $^{2,1,3}$, Dominik Bongartz $^{4,*}$\\

    $^1$ Process Systems Engineering (AVT.SVT), RWTH Aachen University, \\52074 Aachen, Germany\\
    $^2$ JARA-ENERGY, 52056 Aachen, Germany\\
    $^3$ Energy Systems Engineering (IEK-10), Forschungszentrum Jülich, \\52425 Jülich, Germany\\
    $^4$ Department of Chemical Engineering, KU Leuven, \\3001 Leuven, Belgium\\
    $^*$ Corresponding author at: Department of Chemical Engineering, KU Leuven, \\
    3001 Leuven, Belgium. E-mail address: dominikbongartz@alum.mit.edu }
\date{May 2023}

\begin{document}
\maketitle
\textbf{Abstract}	\\
The synthesis of methanol from captured carbon dioxide and green hydrogen could be a promising replacement for the current fossil-based production.
The major energy input and cost driver for such a process is the electricity for hydrogen production. Time-variable electricity cost or availability thus motivates flexible operation. However, it is unclear if each unit of the process should be operated flexibly, and if storage of electricity or hydrogen reduces the methanol production cost.
To answer these questions, we modeled a Power-to-Methanol plant with batteries and hydrogen storage. Using this model, we solved a combined design and scheduling optimization problem, which provides the optimal size of the units of the plant and their optimal (quasi-stationary) operation.
The annualized cost of methanol was minimized for a grid-connected and a stand-alone case study.
The optimization results confirm that storage, especially hydrogen storage, is particularly beneficial when the electricity price is high and highly fluctuating. 
Irrespective of the presence of storage, the whole Power-to-Methanol plant should be operated flexibly: even moderate flexibility of the methanol synthesis unit significantly reduces the production cost. \\

\textbf{Keywords:}
Power-to-Methanol; Flexibility; Hydrogen storage; Battery; Combined design and scheduling optimization; Specific cost of methanol.\\
\\
\textbf{Highlights:}
\begin{enumerate}
\item Simultaneously optimized design and scheduling for Power-to-Methanol 
\item Analyzed the role of flexibility and storage on methanol production cost
\item Already moderate flexibility of the methanol unit significantly reduces methanol cost
\item Storage reduces methanol cost with high and highly variable electricity cost profiles
\item Power-to-Methanol plants should be operated flexibly even with intermediate storage 
\end{enumerate}
\newpage
\section{Introduction}
Methanol is an important platform chemical that is used for the synthesis of several products, e.g., formaldehyde, dimethyl ether, and plastics and has  high potential as alternative fuel or for the production of other liquid fuels, e.g, gasoline and jet fuel \cite{IRENA2021}. 

Nowadays, methanol production mostly relies on fossil feedstocks, e.g., natural gas and coal \cite{IRENA2021}. Methanol from biomass (bio-methanol) and electricity-based methanol (e-methanol) could help achieve environmental objectives. Bio-methanol production is, however, limited by the availability of residual biomass, which could also be used to produce more complex molecules and products \cite{Ulonska2018}\cite{Khadraoui2022b}. E-methanol production relies on a \ce{CO2} source and renewable electricity availability: \ce{CO2} can be obtained from point sources or from air  \cite{Assen2016}, while renewable electricity can be supplied by stand-alone and grid-connected generation parks. Because of the expected increase of installed renewables \cite{InternationalEnergyAgencyIEA2022} and the abundance of \ce{CO2} sources, e-methanol production (Power-to-Methanol) is considered in this work.

Power-to-Methanol plants require high electricity input to produce hydrogen via electrolysis. Because of this high electricity demand, demand-side management \cite{Zhang2016} could be applied to reduce the production cost or to synthesize the product according to the power availability \cite{Burre2020}. However, time-variable operation requires oversizing the plant for a fixed overall production rate. Thus, the extent to which the plant is operated flexibly has to be chosen carefully to balance operating and capital costs \cite{Mucci2023}. 
Oversizing of the electricity-demanding units combined with a product or intermediate storage has often been proposed to deal with renewables or fluctuating prices in different applications, e.g., chlor-alkali electrolysis \cite{Roh2019}, Power-to-Fuel plants \cite{Gorre2019}\cite{Tremel2018}, and ammonia production processes \cite{Allman2018}\cite{Wang2020}. Storage technologies, e.g., batteries and tanks for intermediates, could therefore contribute to reducing the production cost of methanol. In fact, batteries can decouple electricity production from its utilization, while intermediate hydrogen storage can decouple hydrogen production from its conversion into methanol.

While batteries and hydrogen storage systems are frequently used and optimized in microgrids \cite{Valverde2016} to cope with renewables, few works about hydrogen storage and plant capacity optimization are found in the Power-to-Chemical field, for example, for Power-to-Gas  \cite{Gorre2020} and Power-to-Ammonia \cite{Li2020a}\cite{SchulteBeerbuhl2015} plants. The interplay between batteries and hydrogen storage has been investigated for Power-to-Ammonia plants: both storage technologies were chosen by the Cplex optimizer for a case study in which the production rate of ammonia was assumed constant \cite{Osman2020}, while either the hydrogen storage or both storage technologies were chosen according to the solving method and case study in another work \cite{Allman2019}. 
As regards Power-to-Methanol, Chen and Yang \cite{Chen2021} optimized the design of a plant with multiple storage technologies, e.g., hydrogen and thermal storage, for several scenarios. Furthermore, a fuel cell and hydrogen turbine system was included to supply electricity to the plant in case of renewable electricity shortage or expensive dispatch power \cite{Chen2021}. 
However, the Power-to-Methanol configuration with both battery and hydrogen storage has not been investigated. Also, Chen and Yang \cite{Chen2021} linearized the model of highly nonlinear components of the Power-to-Methanol plant with storage and the objective function, thus making, on the one hand, the optimization problem easier to solve, on the other hand, the results less accurate.  
Therefore, while the economic benefits of demand-side management and hydrogen storage utilization in Power-to-Methanol plants are well-known \cite{Tremel2018}\cite{Chen2021}, it is not clear whether batteries and hydrogen storage are always beneficial from an economic perspective, especially in Power-to-Methanol plants, and if and to what extent the methanol synthesis unit should be operated flexibly even in their presence. 

\begin{figure}[h!bt]
    \centering
    \includegraphics[height=0.42\textwidth]{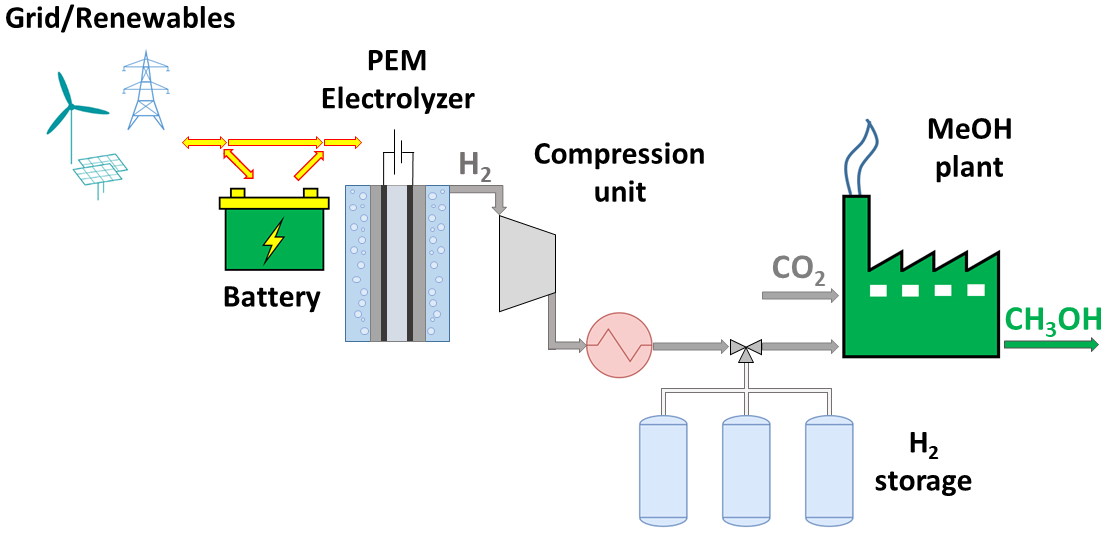}    
    \caption{Sketch of the Power-to-Methanol plant with battery and hydrogen storage. The plant configurations with only the battery, with only the hydrogen storage, and without any storage were also considered (not shown here).}
    \label{fig:Sketch_MEOH_B_V}
\end{figure}

To tackle these questions, we modeled a Power-to-Methanol plant with both a battery and hydrogen storage (Fig. \ref{fig:Sketch_MEOH_B_V}) in GAMS \cite{GAMS} based on mass and energy balances and performance models for the units. The resulting mixed-integer nonlinear programming (MINLP) problem considers combined design (in particular sizing) and scheduling, accounting for the operating limits of every unit, thus fully exploiting the potential of flexible operation \cite{Mitsos2018}.  
We then optimized with BARON \cite{BARON_2018} the design for single scheduling scenarios to investigate how the optimal plant design depends on these. 

The remainder of the paper is structured as follows: in Sections \ref{System} and \ref{Model}, the investigated Power-to-Methanol plant and its model are described; Sections \ref{Case_studies} and \ref{Optimization} present the case studies and the formulation of the optimization problems, while in Section \ref{Results} and \ref{Conclusions}, the optimization results are analyzed and the conclusions are drawn. Furthermore, additional information and the GAMS code are provided in Appendix \ref{Appendix} and Supplementary Information, respectively.
\section{System description} \label{System}
The investigated Power-to-Methanol concept (Fig. \ref{fig:Sketch_MEOH_B_V}) is composed of three key units, i.e., the water electrolyzer, the compression unit, and the methanol synthesis plant. Additionally, a battery and hydrogen storage are considered to investigate their role.
A short description of these units is provided in the following subsections.
\subsection{Battery}
Electricity can be stored electrochemically in batteries in order to be used in following periods. Different battery technologies have been developed for energy storage and automotive purposes \cite{IRENA2017}. For the current work, we considered a lithium-ion type battery because of the high roundtrip efficiency (around 95 \%) and the expected cost drop \cite{IRENA2017}\cite{IEA2021}. 
\subsection{PEM water electrolyzer}
The water electrolyzer is the key component of the investigated Power-to-Methanol concept. Among alternative electrolysis technologies, we chose Polymer Electrolyte Membrane water electrolyzers (PEM-WE) since they have high technological maturity and efficiency, and they are reported to be particularly suitable to follow variable power inputs \cite{Carmo2013}\cite{Kopp2017}. 

Furthermore, PEM-WEs can be operated at relatively high pressures. On the one hand, pressurized electrolyzers have a lower efficiency (higher operating costs); on the other hand, they reduce the size and cost of the downstream compression unit. The trade-off between higher capital or operating cost depends on the considered boundary conditions, e.g., electricity cost and delivery hydrogen pressure, and will be investigated for the considered Power-to-Methanol plant.
\subsection{Hydrogen compression unit}
Hydrogen has to be compressed from the pressure at the outlet of the electrolyzer up to the delivery pressure. Though several compression technologies are available, e.g., mechanical and electrochemical compression, only reciprocating and centrifugal hydrogen compressors have reached a high technological maturity on a large scale \cite{Tahan2022}.
We selected centrifugal compressors because of their high reliability and capability of handling high hydrogen flow rates.

The delivery pressure of hydrogen that the compression unit has to reach varies according to the plant configuration. In absence of hydrogen storage, the delivery pressure matches the operating pressure of the methanol synthesis plant (75 bar). In case intermediate \ce{H2} storage is present, the delivery pressure varies over time according to the pressure inside the storage. 

For the considered case studies, a compression unit without intermediate cooling is sufficient to supply hydrogen at the delivery pressure. The compressed hydrogen is then cooled to 25 °C before being injected into the hydrogen storage or supplied to the methanol synthesis plant.
\subsection{Hydrogen storage}
Hydrogen can be stored either physically, i.e., as liquid or gas, or in some materials via chemical or physical sorption \cite{Moradi2019}. We chose physical-based storage for the considered application since it is the most mature technology \cite{Moradi2019}.

Storing hydrogen as a liquid allows for increasing its energy density, thus reducing the size of the storage. However, this solution requires a dedicated liquefaction plant, which significantly increases the overall plant complexity and cost. 
Instead, gaseous storage does not need dedicated plants except for a compression unit. There are two main gaseous storage alternatives, i.e., above-ground in vessels and underground in caverns \cite{Moradi2019}\cite{Elberry2021}. The use of pressurized vessels was chosen since it is geographically unconstrained. Among the several types of above-ground hydrogen storage \cite{Moradi2019}\cite{Hirscher2010}, hydrogen vessels of type I \cite{Moradi2019} are the most suitable for stationary applications for the pressure range of interest (maximum pressure below 200 bar) as they are the cheapest.
\subsection{Methanol synthesis plant}  
Methanol can be produced from carbon dioxide and hydrogen either directly or indirectly. In the direct route, \ce{CO2} is directly converted to methanol via hydrogenation. In the indirect route, \ce{CO2} is first reduced to carbon monoxide either thermochemically or electrochemically and then mixed with hydrogen to produce syngas that is converted to methanol. The intermediate production of syngas aims at using the consolidated know-how of methanol production from fossil fuels, e.g., natural gas. 
The direct pathway is gaining relevance in the last decades because of the lower plant complexity, overall efficiency, and economic feasibility  \cite{Anicic2014}. For this reason, the direct hydrogenation of \ce{CO2} to methanol pathway was chosen.
\section{Model} \label{Model}
The Power-to-Methanol plant with storage is modeled via mass and energy balances and by considering discrete time dynamics. We assumed quasi-stationary operation for each unit except for storage. Also, cost functions are used to estimate capital and operating costs.
In this section, the main modeling equations of each unit are presented. For the complete model, we refer to the GAMS code provided in the Supplementary Information.

\subsection{Battery}
The battery is modeled via the energy balance 

\begin{equation*} \mathrm{E_{b}}(t) = \mathrm{E_{b}}(t-\Delta t) \cdot (1-r_{\mathrm{self-disch}}) + P_{\mathrm{in}}(t)  \cdot \eta_{\mathrm{ch}} \cdot \Delta t - \frac{P_{\mathrm{out}}(t)}{\eta_{\mathrm{disch}}} \cdot \Delta t, \end{equation*}

where $\mathrm{E_{b}}(t)$ is the energy content, $P_{\mathrm{in}}(t)$ the power input from the grid, $P_{\mathrm{out}}(t)$ the power in DC that can be effectively used, and $\Delta t$ the discretization time step.
The values for the charging ($\eta_{\mathrm{ch}}$) and discharging ($\eta_{\mathrm{disch}}$) efficiencies and the self-discharge rate ($r_{\mathrm{self-disch}}$) were considered constant (see Supplementary Information). The hourly self-discharge rate was estimated by considering an energy loss of 0.2 \% per day \cite{IRENA2017}.
Other modeling approaches, e.g., efficiencies as a function of the current and of the state of charge of the battery \cite{Gonzalez-Castellanos2020}, would have improved the accuracy but also increased the complexity of the optimization problem. 

The power output of the battery can be supplied to the water electrolyzer ($P_{\mathrm{to-PEM}}(t)$) or to the grid ($P_{\mathrm{to-grid}}(t)$) if converted from DC to AC:

\begin{equation*}P_{\mathrm{out}}(t) = \frac{P_{\mathrm{to-grid}}(t)}{\eta_{\mathrm{DC-AC}}} + P_{\mathrm{to-PEM}}(t). \end{equation*}

The capital cost of the battery was estimated from its nominal capacity. The specific cost of 310 \$/kWh (the average installation cost of utility-scale stationary battery systems in 2020  \cite{IEA2021}) was considered and updated to the reference year (2021).
\subsection{PEM water electrolyzer} \label{PEMWE}
The water electrolysis unit is composed of several electrolysis modules. The main specifications of the considered electrolyzer module are summarized in Table \ref{tab:PEMWE_SE_spec}.
The PEM water electrolyzer module was modeled by considering the equations and parameters in Järvinen et al. \cite{Jaervinen2022}. Two additional assumptions were made: i) the Faradaic efficiency is equal to one, and ii) the pressure dependency is considered in the Nernst equation only.

Since the equation-based model contains several variables that are unnecessary for the current scope and would slow down the optimizer, we replaced the electrolyzer model by a polynomial fit of its efficiency (based on the lower heating value, see Appendix \ref{PEMWE_efficiency}). The electrolyzer efficiency was then used to calculate the hydrogen mass flow rate, the power, and cooling demand.

Additional operating costs were considered for the auxiliaries, e.g., the power electronics units and pumps. In particular, the additional energy consumption was assumed equal to the 5\% of the electricity needed to run the electrolyzer.

The capital cost of the PEM electrolyzer was calculated by multiplying the nominal power of the electrolyzer by the specific cost in \texteuro/kW, which was calculated with the correlation proposed by Reksten et al. \cite{Reksten2022}. To also take into account civil works and installation costs, a multiplying factor of 1.75 was assumed. Additionally, a second electrolysis unit is considered to be purchased in the tenth year, since the typical lifetime of PEM electrolyzers (30-90 kh \cite{IEA2019}) is shorter than the typical lifetime of the chemical plant (20-30 years).

\begin{table}[h!bt]
    \centering
    \caption{Main specifications of the considered PEM water electrolyzer module.}
    \begin{tabular}{lc}
        \hline
         & \textbf{Values} \\
          \hline
        Operating temperature &  75 °C \\
        Operating pressure range & 20-40 bar \\
        Nominal power per module & 2.0 MW\\
        Maximum power per module & 2.5 MW\\
         \hline
    \end{tabular}    
    \label{tab:PEMWE_SE_spec}
\end{table}
\subsection{Hydrogen compression unit} \label{H2_compression}
A multi-stage centrifugal compressor was considered in the plant. This compression technology has an operating range of around 50\%-100\% \cite{Tahan2022} if variable rotating speed \cite{He2018} or variable inlet guide vanes \cite{Frank2022} are used. As the other units have a wider operating range, the operation of two compressors in parallel was considered to extend it to 25\%-100\%, although this affects the capital cost of the unit (economy of scale). The use of a low-pressure buffer tank before the compression unit was not considered, although it would have further increased  the flexibility of the plant.

The overall power of compression ($ P_{\mathrm{comp}}$) is calculated as follows:

\begin{equation*}P_{\mathrm{comp}}(t)=\dot{m}_{\mathrm{H_{2}}}(t)\cdot \frac{c_{\mathrm{p}} \cdot T_{\mathrm{in}}}{\eta_{\mathrm{is}} \cdot \eta_{\mathrm{mec}}} \cdot\left(\beta^{\left(\frac{k-1}{k}\right)}-1\right) , \end{equation*}

where $\dot{m}_{\mathrm{H_{2}}}(t)$, $c_{\mathrm{p}} $, and $k$ are the hydrogen mass flow rate, the specific heat at constant pressure, and the heat capacity ratio of hydrogen, respectively, $T_{\mathrm{in}}$ the hydrogen temperature at the inlet of the compression unit (assumed equal to 298 K), $\beta$  the pressure ratio, and $\eta_{\mathrm{is}}$ and $\eta_{\mathrm{mec}}$ the isentropic and mechanical efficiencies. The isentropic efficiency was assumed constant over the operating range and conservatively equal to 0.8 \cite{Tahan2022}\cite{He2018}. 
A validation of the model of the compressor can be found in Appendix \ref{Properties}.

The pressure ratio is determined by the operating pressure of the electrolyzer and the pressure inside the hydrogen storage unless a bypass of the hydrogen storage is considered. With a bypass, some compression energy could be saved depending on the electricity price, the operating load of the methanol synthesis plant, and the pressure inside the storage. However, this configuration was not considered in this work. 
Also, no hydrogen losses were considered during the compression phase since they are lower than 0.5 \% \cite{Tahan2022}.

The compressed hydrogen is then cooled in a single heat exchanger to 298 K. The cooling demand was assumed equal to the power needed for compression.
The capital cost of these components was estimated at the maximum flow rate condition by using the cost models proposed by Biegler et al. \cite{Biegler1997}. 
\subsection{Hydrogen storage}
The amount of hydrogen available in the vessel ($M_{\mathrm{H_{2}}}(t)$) is calculated with the following mass balance:

\begin{equation*}M_{\mathrm{H_{2}}}(t) = M_{\mathrm{H_{2}}}(t- \Delta t) + \dot{m}_{\mathrm{H_{2,prod}}}(t) \cdot \Delta t - \dot{m}_{\mathrm{H_{2,MeOH}}}(t) \cdot \Delta t, \end{equation*}

where $\dot{m}_{\mathrm{H_{2,prod}}}(t)$ is the hydrogen flow rate produced by the electrolyzer,  and $\dot{m}_{\mathrm{H_{2,MeOH}}}(t)$ the hydrogen flow rate consumed in the methanol synthesis plant.

The hydrogen storage was assumed isothermal at 298 K, and the pressure inside the storage ($p(t)$) was estimated from the available hydrogen inside the storage by using a linear interpolation of the density ($\rho$) (Appendix \ref{Properties}).
As the methanol synthesis occurs at a pressure of 75 bar, the minimum pressure of the storage is set to that value. The available hydrogen stored in the vessel at the time step \textit{t} and its pressure are, therefore, calculated as follows:

\begin{equation*} M_{\mathrm{H_{2}}}(t)= V\cdot\left(\rho_{\mathrm{(p(t))}}-\rho_{\mathrm{(75\; bar)}}\right) \approx V \cdot 0.073\;\si{ \frac{kg}{m^3 bar}} \cdot \left(p(t)-75\; \si{bar}\right), \end{equation*}  

where \textit{V} is the volume of the storage, which was considered as a continuous variable, even though multiple hydrogen vessels of limited size would be probably needed in large-scale plants \cite{Elberry2021}.

The capital cost of the storage was estimated by using the specific cost per unit of volume. Further details can be found in Appendix \ref{CAPEX_vessel}.
\subsection{Methanol synthesis plant} 
The methanol synthesis plant was modeled in Aspen Plus \textsuperscript{\textregistered} V11 (see Appendix \ref{Methanol} for more detail). However, such a detailed model cannot be embedded in GAMS and handled by the optimizer. Therefore, the key information about efficiency, mass flow rates, and energy demand was extracted to build a simplified algebraic model of the plant.
In particular, the produced methanol mass flow rate depends linearly on the hydrogen flow rate, while the electricity demand, the cooling demand, and the inlet \ce{CO2} mass flow rate scale linearly with the methanol production rate. 
The complete model can be found in the GAMS file provided in the Supplementary Information.

While the \ce{H2} production is simulated by the model of the electrolysis unit (Section \ref{PEMWE}), the model of the carbon capture unit for the \ce{CO2} supply was not included in the optimization problem. Therefore, it was assumed that a suitable \ce{CO2} stream can be supplied to the methanol synthesis plant at every time step.

The capital cost of the methanol synthesis plant ($\mathrm{CAPEX}$) was estimated by using the Aspen Plus model (Appendix \ref{Methanol}) as the reference in the regression function, which is defined as:

\begin{equation*}\mathrm{CAPEX= CAPEX_{ref}\cdot \left(\frac{M_{MeOH}}{M_{MeOH,\;ref}}\right)^b} , \end{equation*}

where $b$ is the regressive exponent (equal to 0.6 for the six-tenths rule), to evaluate the $\mathrm{CAPEX}$ for different plant capacities ($\mathrm{M_{MeOH}}$). Further details about the cost estimation can be found in Appendix \ref{CAPEX_MeOH}.
\section{Case studies} \label{Case_studies}
To investigate the effect of storage and flexibility on the Power-to-Methanol plant, two case studies were considered, i.e., grid-connected and stand-alone plants.

For both case studies, two-month scenarios were considered and discretized with hourly resolution (1440 time steps) for the scheduling problem. The relatively short time frame is motivated by the size of the optimization problems (the number of variables and constraints scale linearly with the number of time steps unless time-aggregation series techniques are used \cite{Schafer2020}) and by the fact that no seasonal storage of hydrogen and electricity is expected. 
Nevertheless, this time frame allows considering price and renewable power variations that occur on different time scales, e.g., electricity price variation between working and weekend days for multiple weeks, and the diurnal effects on renewable power generation. 
Furthermore, as the plant units have high ramp rates and fast dynamics \cite{IRENA2017}\cite{Kopp2017}\cite{Dieterich2020}\cite{FlexMeOH22} compared with the considered time step (1 h), quasi-stationarity at the operating points was assumed: the plant is considered in steady-state at each time step independently from the operating point at the previous time step.

\begin{figure}[h!bt]
    \centering
    \includegraphics[width=1\textwidth]{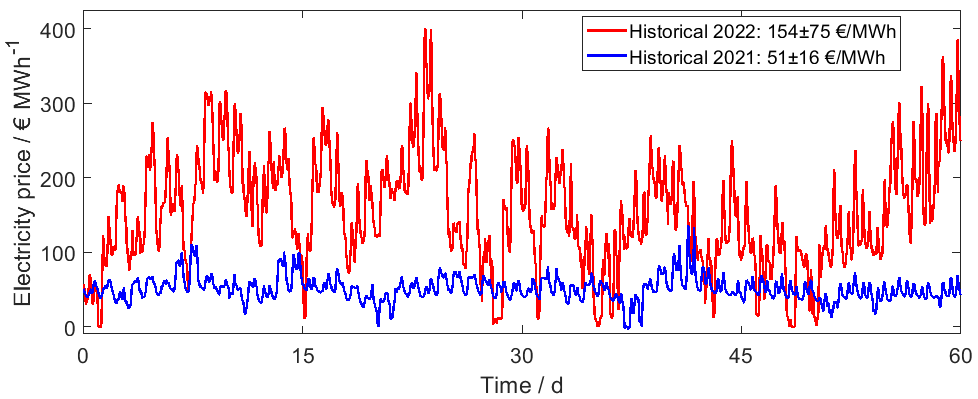}    
    \caption{Electricity price scenarios \cite{ENTSOE}.}
    \label{fig:Electricity_price_scenarios}
\end{figure}

For the grid-connected case study, we assumed that there are no limitations in terms of power availability, and that the amount of energy exchanged with the grid does not affect the electricity price. Two electricity price scenarios based on historical day-ahead electricity prices of Germany (the first 60 days of the years 2021 and 2022 \cite{ENTSOE}) were analyzed. 
The scenarios differ in terms of mean value and standard deviation (Fig. \ref{fig:Electricity_price_scenarios}), representing a low-price scenario with a small standard deviation, and a high-price scenario with a high standard deviation.

\begin{figure}[h!bt]
    \centering
    \includegraphics[width=1\textwidth]{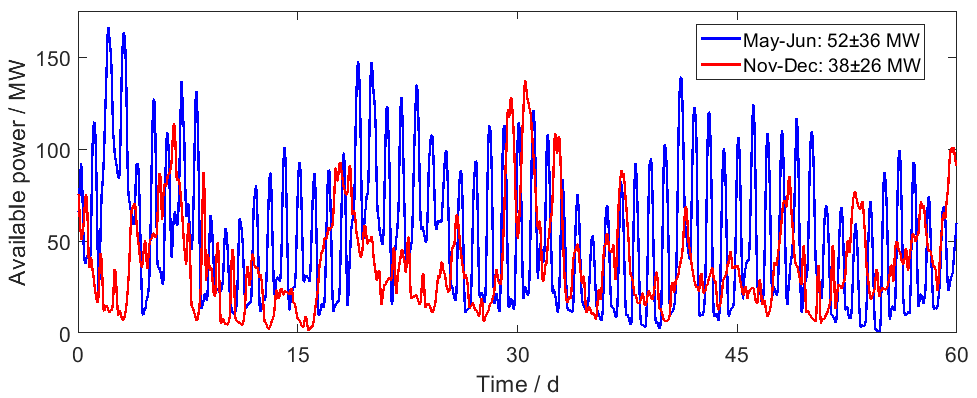}    
    \caption{Power availability scenarios \cite{ENTSOE}.}
    \label{fig:Power_availability}
\end{figure}

For the stand-alone case study, the power is produced from a renewable generation park composed of wind and photovoltaic plants. The power generation profiles (Fig. \ref{fig:Power_availability}) were taken from German historical data (May-June 2022 and November-December 2022 \cite{ENTSOE}) and scaled to a typical size for large renewable parks (126 MW and 110 MW in nominal conditions for wind and photovoltaic plants, respectively). 
Additional assumptions were made: i) the forecast of the energy production from the renewable generation park is perfect, ii) the cost of electricity is null since the capital cost of the renewable plants \cite{IRENA_RenGen_2021} is included in the initial investment cost, iii) the sum of the electricity demand of the whole Power-to-Methanol plant must be less than or equal to the electricity available for each time step; iv) electricity in the battery can also be converted to AC to supply energy to the methanol synthesis plant in case of renewable production shortage; v) the eventual excess of electricity is curtailed.

\begin{figure}[h!bt]
    \centering
    \includegraphics[width=1\textwidth]{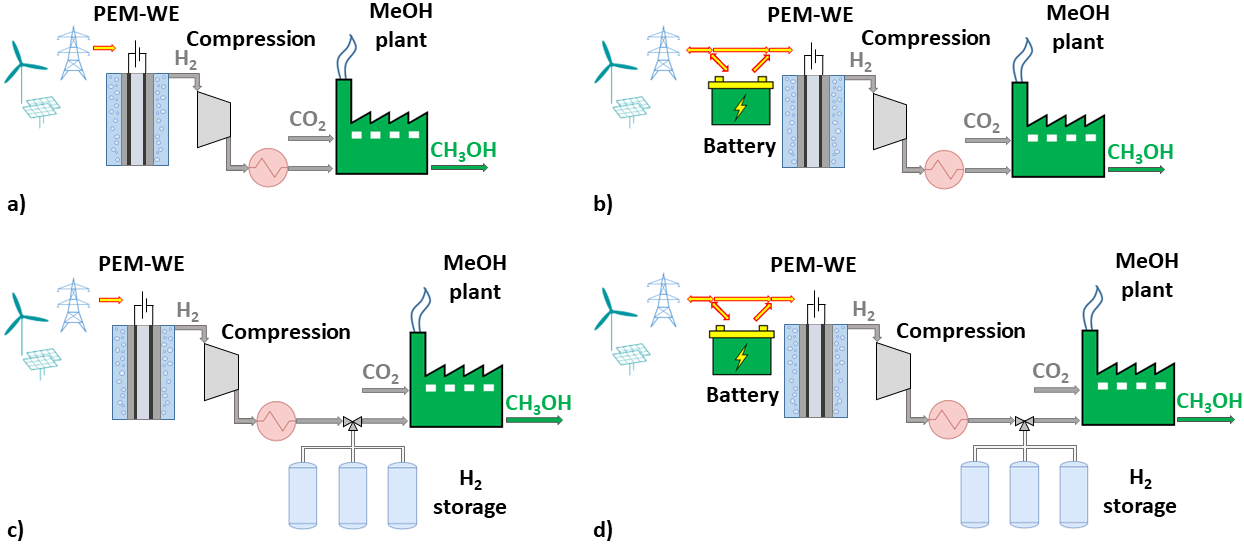}    
    \caption{Investigated Power-to-Methanol plant configurations: a) without any storage, b) with a battery, c) with \ce{H2} storage, and d) with battery and \ce{H2} storage.}
    \label{fig:4_options}
\end{figure}

To investigate the role of storage, four plant configurations were considered, i.e., the Power-to-Methanol plant without storage, with a battery, with \ce{H2} storage, and with both storage technologies (Figure \ref{fig:4_options}) in the grid-connected case study. In the stand-alone case study instead, only the Power-to-Methanol plant configuration with both storage technologies was considered since storage is essential to grant continuous operation of the downstream plant. In fact, the battery can supply energy to the plant in case of a temporary electricity generation shortage, while the \ce{H2} vessel can supply \ce{H2} to the methanol synthesis plant when the electricity production is not sufficient to run the electrolysis unit.

Finally, the annualized cost of methanol ($C_{\mathrm{MeOH}}$) was used as a metric to compare the configurations for each case study. It is defined as:

\begin{equation*}C_{\mathrm{MeOH}}=\mathrm{\frac{CAPEX_0+CAPEX_{\, PEM} \cdot (1+\textit{i})^{-10} + \sum_{j=1}^N ( O\&M+OPEX_{\, y})\cdot (1+\textit{i})^{\textit{-j}}}{\sum_{j=1}^N ( \textit{M}_{MeOH,\; y})\cdot (1+\textit{i})^{\textit{-j}}} }, \end{equation*}

where $\mathrm{CAPEX_0}$, O\&M, and $\mathrm{OPEX_{\, y}}$ are the initial capital expenditures, the yearly operation and maintenance, and the operating costs of the plant, respectively, N the number of operating years, $M_{\mathrm{MeOH,\; y}}$ the yearly production of methanol, and $i$ the discount rate. Additionally, a replacement of the electrolysis unit at the end of the tenth year is considered (see Section \ref{PEMWE}).
The capital cost of the single units was calculated with the cost models described in Section \ref{Model} as the result of the optimal design, while the operating costs and the methanol production as the result of the optimal scheduling. Further assumptions for the economic evaluations are summarized in Table \ref{tab:Assumptions}. 

\begin{table}[h!bt]
    \centering
    \caption{Main assumptions for the economic evaluations.}
    \begin{tabular}{lc}
        \hline
         & \textbf{Values} \\
          \hline
        Operating weeks per year for the methanol synthesis plant &  48 w/y\\
        Lifetime of the Power-to-Methanol plant &  20 y\\
        Yearly O\&M as \% of the initial $\mathrm{CAPEX_0}$ & 5 \%\\
        Discount rate $i$ & 5 \%\\
        \$ to \texteuro \;  conversion factor & 0.85\\
        Cooling cost & 1  \texteuro/MWh\\
        \ce{CO2} stream cost & 50  \texteuro/t \cite{IRENA2021}\\
         \hline
    \end{tabular}    
    \label{tab:Assumptions}
\end{table}
\section{Optimization} \label{Optimization}
The combined design and scheduling optimization problems were formulated in GAMS 38.2.1 \cite{GAMS}.  
These MINLP problems were initialized with a starting solution found after a multi-start local search and then solved with the global optimizer BARON, version 22.2.3 \cite{BARON_2018}. Two stopping criteria were considered for the optimizer, i.e., an optimality gap of $0.01$ and the CPU time limit of 12 h. Unfortunately, the relative difference between the upper bound (best feasible solution found) and the lower bound (best possible solution) did not reach the required optimality gap within the computational time limit in most of the optimization problems, thus, overall global optimality cannot be guaranteed. Nevertheless, the upper bound stabilized after a few iterations, thus suggesting that the convergence issues were mainly due to the difficulty of proving global optimality with wide variable bounds (global optima were proven in some of the problems after tightening the variable bounds). In preliminary analyses, we also tested Knitro, ANTIGONE, SCIP, and DICOPT, but their performance was worse than that of BARON. 

For all the considered case studies, the annualized cost of methanol ($C_{\mathrm{MeOH}}$) was minimized. 
In the following subsections, the optimization variables and the main constraints are discussed.
\subsection{Optimization variables} \label{Optimization_variables}
The considered optimization variables can be divided into two categories: design and scheduling variables. The design variables allow defining the optimal size and capacities of the units, which mainly contribute to the investment cost, while the scheduling variables how the units are operated, thus affecting the operating costs and the production rates. 

 As regards the plant design, the size of the battery, electrolyzer, compression unit, \ce{H2} storage, and methanol synthesis plant was optimized. In particular, the size of the methanol synthesis plant was estimated from the nominal \ce{H2} input stream that can differ from the maximum hydrogen production rate of the electrolysis unit in case there is \ce{H2} storage.
Two additional design variables that also have a high impact on the operating phase were considered, i.e., the operating pressure of the electrolysis unit and the maximum pressure ratio of the compressor. The first affects the electrolyzer efficiency as well as the size of the downstream compression unit, while the second determines the maximum pressure that can be achieved in the \ce{H2} storage. The choice of these two design variables aims at capturing the interplay between the size of the compression unit and the \ce{H2} storage.

A summary of the design variables and their bounds can be found in Table \ref{tab:Opt_vars}.  
A maximum size of 100 MW (40 modules) was considered as representative of the electrolyzer size in the short-term future \cite{Reksten2022}. The maximum size of the battery and \ce{H2} storage was defined accordingly in order to provide a few hours of storage (around 4 h and 12 h, respectively, when the electrolyzer size is 100 MW). The electrolyzer pressure range is typical for PEM water electrolyzers. The maximum pressure ratio of the compressor unit is limited by technical issues, e.g., discharge temperature \cite{Tahan2022} and overall efficiency.

In order to calculate the optimal scheduling, optimization variables for the battery, the electrolyzer, and the \ce{H2} storage were introduced.
For the battery, the  power input ($P_{\mathrm{in}}(t)$) and the  power output ($P_{\mathrm{to-grid}}(t)$ and $P_{\mathrm{to-PEM}}(t)$) were optimized. Additionally, two sets of binary variables ($x_{\mathrm{b, \;in}}(t)$ and $x_{\mathrm{b, \;out}}(t)$) were introduced to define operational constraints. Similarly, the power taken from the grid ($P_{\mathrm{grid}}(t)$) and a binary variable ($x_{\mathrm{on-off}}(t)$) were considered for the electrolysis unit. 
The hydrogen flow rate fed into the methanol synthesis plant ($\dot{m}_{\mathrm{H_{2,MeOH}}}(t)$) determines the methanol production rate over time. 
The bounds of the operational variables for optimization depend on both the design optimization variables and technical constraints and will be discussed in more detail in Section \ref{Operation_constraints}. 

The variables in Table \ref{tab:Opt_vars} are shown for the more general configuration of the Power-to-Methanol plant with both storage technologies. Though the four investigated configurations (Figure \ref{fig:4_options}) could have been optimized via a single superstructure optimization problem by introducing binary variables, enumeration was preferred in order to limit the complexity of the optimization problem, which already counts 11526 variables (6 design and 11520 scheduling variables) for the flexible Power-to-Methanol plant with both storage technologies.
The optimization problem was modified when considering simpler plant configurations (no storage and no flexibility) to reduce the problem size and improve convergence. 
\begin{table}[h!bt]
    \centering
    \caption{Design and scheduling variables of the optimization problem for the grid-connected and stand-alone case studies. A short description of the variable and the considered upper and lower bounds are also provided.}
\rotatebox{90}{%
    \begin{tabular}{lcccc}
        \hline
         \textbf{Unit}&\textbf{Symbol}& \textbf{Description of the design variable} & \textbf{Min value}& \textbf{Max value} \\
          \hline
        Battery & $\mathrm{E_{b,\, nom}}$ &  Energy capacity & 5 MWh& 400 MWh \\
        PEM-WE & $N_{\mathrm{mod}}$ &  Number of modules & 10 & 40 \\
                               & $p_{\mathrm{PEM}}$ &  Operating pressure & 20 bar& 40 bar \\
        \ce{H2} compression & $\beta_{\mathrm{max}}$ &  Maximum pressure ratio  & 1.88 & 3.50 \\
        \ce{H2} storage & $V$ &  Volume of the vessel & 25 \si{m^{3}} & 5000 \si{m^{3}}  \\
                        & $\dot{m}_{\mathrm{H_{2,MeOH},nom}}$ &  Nominal \ce{H2} flow rate for the methanol synthesis plant & 0.4 t/h & 2 t/h \\
        \hline
         & & & &  \\
        \hline
         \textbf{Unit}&\textbf{Symbol}& \textbf{Description of the scheduling variable} & \textbf{Min value}& \textbf{Max value} \\
          \hline
        Battery & $P_{\mathrm{in}}(t)$ &  Power input & 0.5 MW * & min(100 MW, $\left(0.9 \cdot \mathrm{E_{b,\, nom}/1h}\right)$) \\
                & $x_{\mathrm{b, \;in}}(t)$ &   Binary to decide when the battery is charged & 0 & 1 \\
                & $x_{\mathrm{b, \;out}}(t)$ &   Binary to decide when the battery is discharged & 0 & 1 \\
                & $P_{\mathrm{to-grid}}(t)$ &  Power output, electricity back to the grid (AC) & 0.5 MW * & min(100 MW, $\left(0.9 \cdot \mathrm{E_{b,\, nom}/1h}\right)$) \\
                & $P_{\mathrm{to-PEM}}(t)$ &  Power output, electricity to the  electrolyzer & 0.5 MW * & min(100 MW, $\left(0.9 \cdot \mathrm{E_{b,\, nom}/1h}\right)$) \\
        PEM-WE & $P_{\mathrm{grid}}(t)$ &  Power directly from the grid/renewable plant & 0.5 MW $\cdot N_{\mathrm{mod}}$ * & 2.5 MW $\cdot N_{\mathrm{mod}}$  \\
                               & $x_{\mathrm{on-off}}(t)$ &  Binary to activate or deactivate the electrolyzer &  0 & 1 \\    
        \ce{H2} storage & $\dot{m}_{\mathrm{H_{2,MeOH}}}(t)$ &  \ce{H2} flow rate for the methanol synthesis plant & $0.2 \cdot \dot{m}_{\mathrm{H_{2,MeOH},nom}}$ & $\dot{m}_{\mathrm{H_{2,MeOH},nom}}$ \\                                       
        \hline
* 0 MW & if not operating\\
        \hline
    \end{tabular}    
 }%
    \label{tab:Opt_vars}
\end{table}
\subsection{Operation constraints} \label{Operation_constraints}
The operation of the Power-to-Methanol plant depends on the size of the units and on their operational limits.

Batteries have wide operating ranges with respect to their state of charge \cite{IRENA2017}. In order to prevent deep discharge phenomena, a minimum state of charge of 10 \% was assumed. Some operational constraints on the hourly energy flows were also set. When active, a minimum power flow of 0.5 MW was considered in order to avoid small energy transfers, which would not be feasible for the power electronics in large plants. The maximum charge and discharge rates were also limited by either the nominal capacity of the battery or the maximum power of the electrolyzer (100 MW).

PEM water electrolyzers also have extremely wide operating limits (0-160 \% relative to the nominal load \cite{IEA2019}). However, we considered an operating range of 25-125 \% since the downstream compression unit has a narrower operating range (25-100 \% of the \ce{H2} flow rate) and no buffer tank was considered between the electrolyzer and the compressor unit. 

The \ce{H2} storage located after the compression unit has a 0-100 \% operating range, where 0\% corresponds to the pressure inside the storage of 75 bar (no \ce{H2} available for the methanol synthesis in the storage). 

The methanol synthesis plant was assumed to always be operating over the considered time horizon within a range of 20-100 \% of the nominal load \cite{Dieterich2020}. This choice also allows avoiding the implementation of additional optimization variables for the on-off, start-up and shut-down phases and the introduction of an economic penalty due to the energy and hydrogen consumption during these phases. 

As mentioned in Section \ref{Optimization_variables}, the methanol flow rate is calculated from the hydrogen flow rate fed into the methanol synthesis plant ($\dot{m}_{\mathrm{H_{2,MeOH}}}(t)$). This hydrogen flow is equal to the hydrogen flow produced by the electrolyzer unless there is the \ce{H2} storage. In this case, the hydrogen flow to the methanol synthesis plant is considered an optimization variable. Furthermore, to investigate the effect of the flexibility of the methanol synthesis plant on its production cost, the following ramp constraint was set:

\begin{equation*} |\dot{m}_{\mathrm{H_{2,MeOH}}}(t)- \dot{m}_{\mathrm{H_{2,MeOH}}}(t-\Delta t)| \leq \mathrm{RL} \cdot \dot{m}_{\mathrm{H_{2,MeOH},nom}}, \end{equation*}

where RL is the maximum ramp limit expressed in percentage, and $\dot{m}_{\mathrm{H_{2,MeOH}},nom}$ is the nominal \ce{H2} flow rate fed into the methanol synthesis plant. In case no flexibility was considered, steady-state operation at the nominal capacity of the methanol synthesis plant was assumed.

A summary of the operating limits and of the operational constraints can be found in Tables \ref{tab:Opt_vars} and \ref{tab:Operation_constr}. We refer to the GAMS files in the Supplementary Information for the implementation of all constraints.

\begin{table}[h!bt]
    \centering
    \caption{Summary of the operation ranges of the plant units (minimum-maximum load).}
\begin{tabular}{lc}
        \hline
         \textbf{Unit}&\textbf{Operating range} \\
          \hline
        Battery & 10-100 \% \\
        PEM-WE & 25-125 \%* \\
        \ce{H2} compression & 25-100 \% \\
        \ce{H2} storage & 0-100 \% \\
        Methanol synthesis plant & 20-100 \% \\
        \hline
        * referred to the nominal load\\
        \hline
\end{tabular} \label{tab:Operation_constr}
\end{table}
\section{Results and discussion} \label{Results}
In this section, the optimization results for the two case studies (grid-connected and stand-alone) are shown. 

\subsection{Grid-connected case study}
In order to investigate the role of storage, we compared the four plant configurations (Fig. \ref{fig:4_options}) for both a low-price low-variable and a high-price highly-variable electricity profile scenarios (Fig. \ref{fig:Electricity_price_scenarios}).
The role of flexibility was investigated by varying the maximum ramp limit (RL) of the methanol synthesis plant (in \%/h of the nominal load).

\subsubsection{Low-price low-variable scenario (Historical 2021)}
\begin{figure}[h!bt]
    \centering
    \includegraphics[width=1\textwidth]{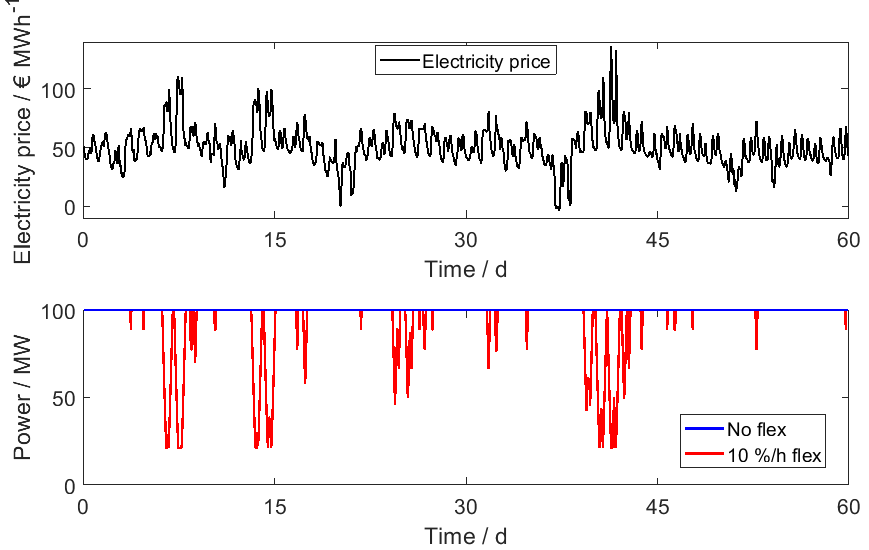}  
    \caption{Power of the PEM electrolyzer for two levels of flexibility of the methanol synthesis plant, i.e., `No flexibility' and `10 \%/h flexibility'. The results are shown for the Power-to-Methanol plant without any storage for the low-price low-variable electricity price scenario (top part of the figure).}
    \label{fig:Power_PEMWE_GC_Scen2021}
\end{figure}

The calculated specific cost for methanol is around 0.9 \texteuro/kg. 
As shown in Table \ref{tab:Optimization_GC_2021}, the specific cost gets slightly lower (around 1 \%) by allowing flexible operation of the methanol synthesis plant. Since the optimal size of the methanol synthesis plant does not change by allowing flexibility in all the configurations for this electricity scenario (the electrolyzer size is at its upper bound), the cost reduction is due to the savings on the operating costs, although the overall methanol production is reduced because of the part-load operation. In fact, if the methanol synthesis plant is operated flexibly, the electrolyzer, the most energy-intense unit, can be operated at the maximum load when the electricity is cheap, and at a lower load when  the electricity is expensive (Fig. \ref{fig:Power_PEMWE_GC_Scen2021}). Nevertheless, the capacity factor of the methanol synthesis plant is above 94 \% (calculated on the working hours) even when the plant is operated flexibly in the considered scenario.

Storage does not improve the economic competitiveness of the plant in this scenario: the optimal solution is without any storage. The potential operating cost savings due to storage do not repay the investment cost for storage because of the relatively low mean value and standard deviation of the electricity price profile.

As regards the other design variables, the optimal pressure of the electrolyzer reaches the upper bound (Table A.\ref{tab:Optimization_GC_2021_complete}), which means that the increased operating costs due to the reduction of the electrolysis efficiency are lower than the additional capital cost that a bigger compression unit would have had.
Also, the nominal hydrogen flow rate for the methanol synthesis plant is equal to the maximum hydrogen production of the electrolyzer since decoupling \ce{H2} storage is not used.  
An optimal design solution close to the upper bounds was expected for not highly-fluctuating electricity price profiles without restrictions on power availability since the cost of the plant components benefits the economy of scale.

\begin{table}[h!bt]
    \centering
    \caption{Optimization results for the grid-connected case study without storage for the low-price low-variable electricity price scenario. Refer to Appendix \ref{Additional_results} for the complete results.}
    \begin{tabular}{ccccc}
        \hline
         \textbf{Flexibility} & \textbf{Specific cost of MeOH}&  \textbf{$\mathrm{MeOH_{y}}$} & \textbf{$\mathrm{CAPEX_{0}}$}& \textbf{$\mathrm{OPEX_{y}}$} \\
          \%/h & \texteuro/kg &   kt & M\texteuro & M\texteuro\\
          \hline
        0 & 0.91    &  80.3  & 133 & 51 \\
        5  & 0.90  &  76.2 & 133  & 46 \\
        10 & 0.90 &  75.7   & 133  & 46 \\
        \hline
       \end{tabular} 
    \label{tab:Optimization_GC_2021}
\end{table}
\FloatBarrier
\subsubsection{High-price high-variable scenario (Historical 2022)}
In this scenario, the methanol cost  is around two times higher than the low-price scenario, thus showing the high influence of the electricity price on the final cost of the product.

\begin{figure}[h!bt]
    \centering
    \includegraphics[height=0.25\textheight]{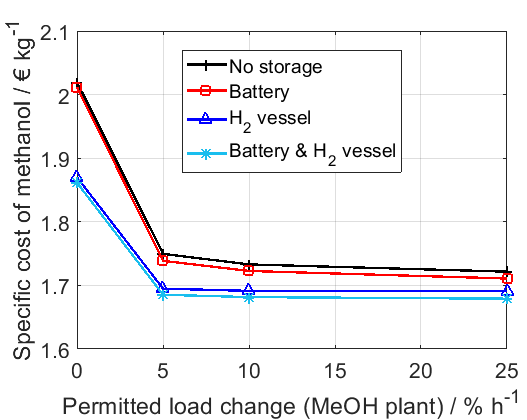}  
    \caption{Specific cost of methanol vs. permitted load change of the methanol synthesis plant for 4 plant configurations, i.e., without any storage, with a battery, with a \ce{H2} vessel, and with both storage for the high-price high-variable electricity price scenario.}
    \label{fig:Cost_vs_flexibility_GC_Scen2022}
\end{figure}

Flexible operation significantly (10-15 \%) reduces the methanol cost compared to the previous scenario. In particular, there is a steep decrease in the methanol cost from the `No flexibility' and `5 \%/h flexibility' cases  (Fig. \ref{fig:Cost_vs_flexibility_GC_Scen2022} and Table \ref{tab:Optimization_GC_2022}). This result is explained by the fact that a relevant part of the hourly price variations is lower than 5\% (see Appendix \ref{Scenarios_Appendix}) and by the possibility of reducing the production rate of methanol during the highly expensive electricity hours that characterize this scenario. This frequent, though optimal, operation at part loads reduces the methanol synthesis plant capacity factor to around 60 \% for all the configurations in this price scenario. Further minor relative reductions of the cost can be noticed when higher ramp limits are allowed since the methanol production rate can  adapt better to the highly fluctuating electricity price profile. 

\begin{figure}[h!bt]
    \centering
    \includegraphics[width=1\textwidth]{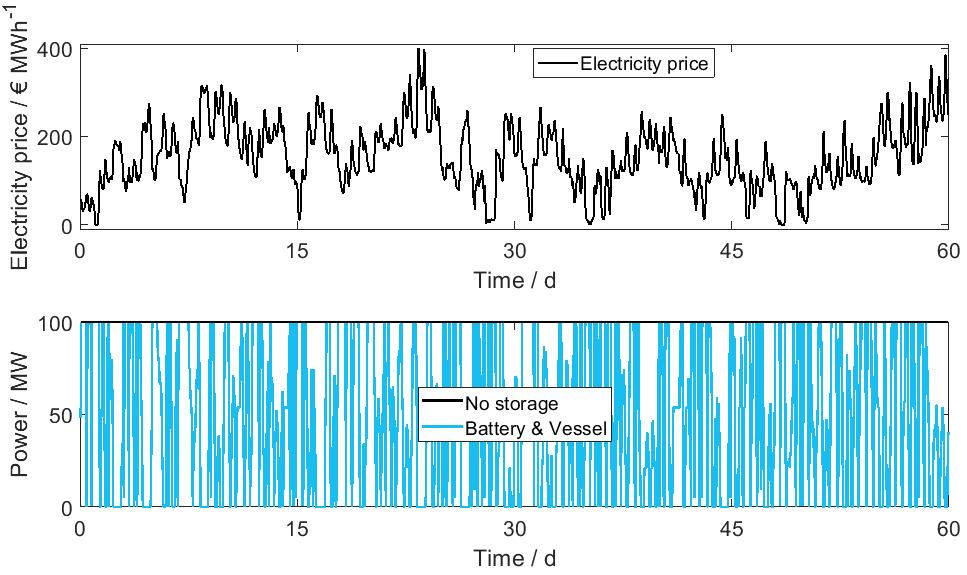}    
    \caption{Power of the PEM-WE taken from the grid over time for the `No flexibility' case. The best (`Battery \& Vessel') and worst (`No storage') cases are shown. The electricity cost profile is also shown to better interpret the optimal scheduling.}
    \label{fig:Power_H2prod_GC_Scen2}
\end{figure}

\begin{figure}[h!bt]
    \centering
    \includegraphics[width=1\textwidth]{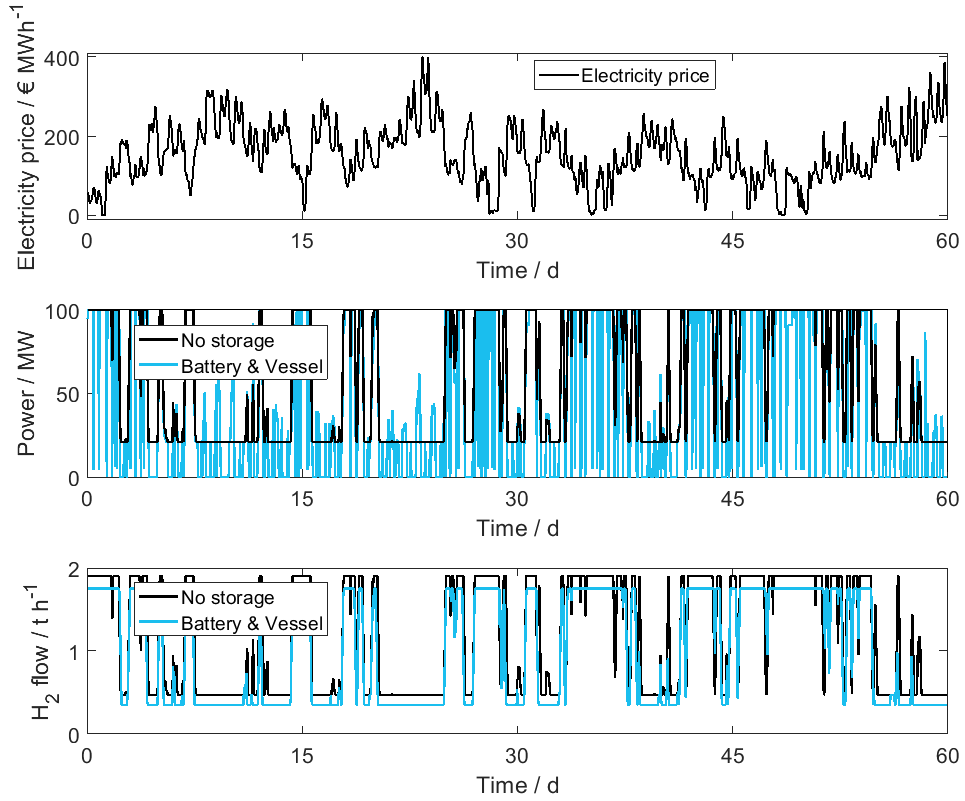}    
    \caption{Power of the PEM-WE taken from the grid (in the middle) and \ce{H2} flow rate to the methanol synthesis plant (bottom part) for the `25 \%/h flexibility' case. The best (`Battery \& Vessel') and worst (`No storage') cases are shown. The electricity cost profile is also shown to better interpret the optimal scheduling.}
    \label{fig:Power_H2used_GC_Scen2}
\end{figure}

In this scenario, storage plays a relevant role, regardless of whether the methanol synthesis plant is operated constantly at nominal load or flexibly. The presence of the hydrogen vessel has the highest impact on the methanol cost since it decouples the production and the consumption of hydrogen and allows downsizing the methanol synthesis plant (7 \% to 40 \% when operating flexibly and constantly at nominal load, respectively). Thanks to the \ce{H2} storage, the electrolyzer can be turned off during expensive electricity hours (Figures \ref{fig:Power_H2prod_GC_Scen2} and \ref{fig:Power_H2used_GC_Scen2}), while the downstream plant is still producing methanol. Nevertheless, the optimal storage size gets lower in case the MeOH plant is operated flexibly, as already observed in literature \cite{Tremel2018}\cite{Li2020a}\cite{Chen2021}. 
Furthermore, the methanol synthesis plant with \ce{H2} storage has a smoother operational profile than the case without \ce{H2} storage (see the bottom part of Fig. \ref{fig:Power_H2used_GC_Scen2}). 
The battery unit also contributes to reducing the methanol cost by both supporting the electrolyzer and selling electricity back to the grid, thus reducing the operating costs.

As regards the design of the Power-to-Methanol plant, the electrolyzer size always reaches the upper bound (see Table A.\ref{tab:Optimization_GC_2022_complete}). The nominal hydrogen flow rate into the methanol plant differs from the maximum amount of hydrogen that can be produced from the electrolyzer (around 1.9 t/h) only when there is \ce{H2} storage.
Interestingly, when there is no \ce{H2} storage and the methanol synthesis plant is operated constantly at nominal load, more efficient production of hydrogen seems to be favored because of the high cost of electricity: The operating pressure is below the upper bound (around 31.5 bar), differently from all the other configurations.


\begin{table}[h!bt]
    \centering
    \caption{Optimization results for the grid-connected case study for the high-price high-variable electricity price scenario. Only the results of the Power-to-Methanol configurations with (B\&V) and without (No B\&V) both storage technologies are shown here. Refer to Appendix \ref{Additional_results} for the complete results.}
    \begin{tabular}{lccccccc}
        \hline
         & \textbf{Flexibility}& \makecell{\textbf{Specific cost} \\ \textbf{of MeOH}}& \textbf{$\mathrm{En_{batt}}$}& \textbf{$V$}  & \textbf{$\mathrm{MeOH_{y}}$} & \textbf{$\mathrm{CAPEX_{0}}$} & \textbf{$\mathrm{OPEX_{y}}$}\\
        & \%/h & \texteuro/kg & MWh & \si{m^{3}} &   kt & M\texteuro & M\texteuro\\
          \hline
        B\&V& 0 &  1.86  & 115  & 4730  &   47.0  & 196 & 57 \\
        No B\&V& 0 & 2.02  & - & - &   80.7  & 135  & 141 \\
        \hline
        B\&V& 5& 1.69  & 116  & 1210  &  42.4  & 179  & 43 \\
        No B\&V& 5& 1.75 & - & - &   49.5 & 133  & 64 \\
        \hline
        B\&V& 25 &1.68 & 111  & 1040  &  42.9 & 177  & 44\\
       No B\&V&25 & 1.72  & - & -  &  49.4  & 133& 63 \\
        \hline
    \end{tabular}    
    \label{tab:Optimization_GC_2022}
\end{table}
\subsection{Stand-alone case study} 
\begin{figure}[h!bt]
    \centering
    \includegraphics[width=1\textwidth]{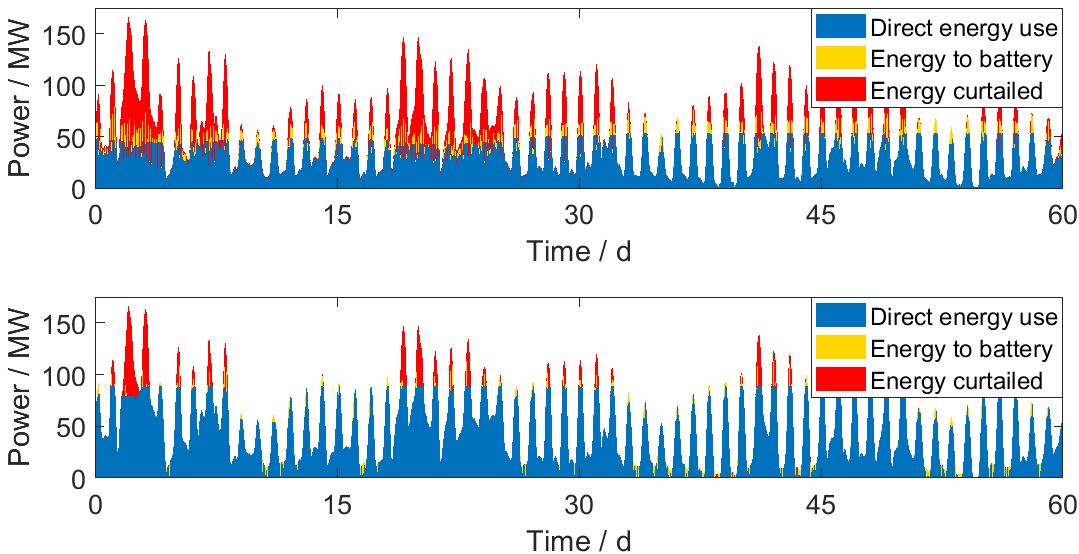}    
    \caption{Energy consumption, energy that is stored in the battery in that hour, and energy curtailed for the `No flexibility' (top) and $\pm$ `25 \%/h flexibility' (bottom) cases (May-Jun scenario).}
    \label{fig:Power_SA_summer}
\end{figure}

\begin{figure}[h!bt]
    \centering
    \includegraphics[width=1\textwidth]{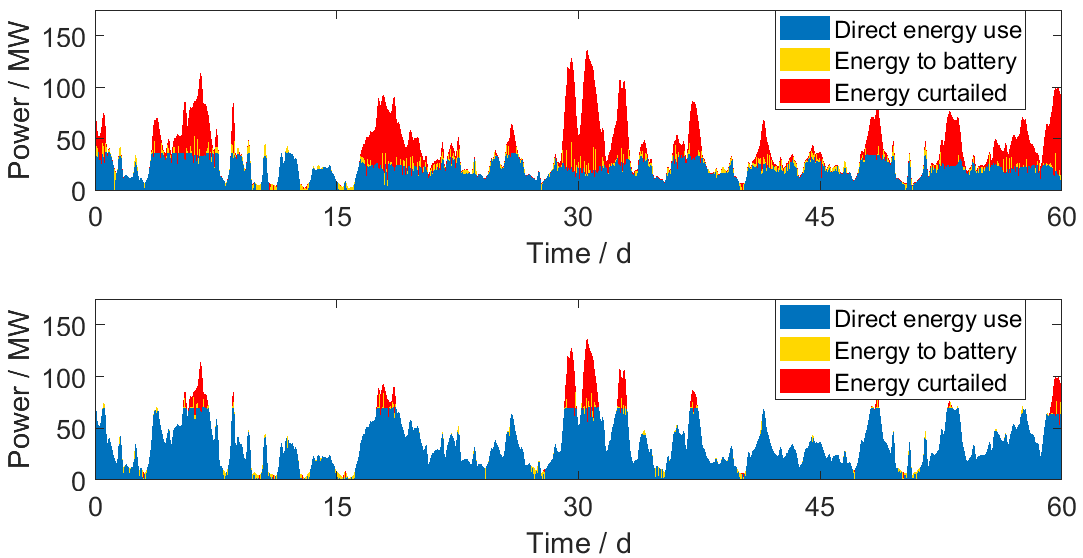}    
    \caption{Energy consumption, energy that is stored in the battery in that hour, and energy curtailed for the `No flexibility' (top) and $\pm$ `25 \%/h flexibility' (bottom) cases (Nov-Dec scenario)}
    \label{fig:Power_SA_winter}
\end{figure}

Also in the stand-alone case study, flexible operation of the methanol synthesis plant highly affects the specific cost of methanol in both scenarios even though both storage technologies are used (Table \ref{tab:Optimization_SA}). The steep reduction of the cost between the `No flexibility' and `5 \% flexibility' case is mainly due to the possibility of increasing the consumption of the available electricity and thus the methanol production, although the capacity factor of the methanol synthesis plant is around 60 \%. In both scenarios, the utilization factor of the produced electricity rises from around 60-70 \% to over 90 \% (Figures \ref{fig:Power_SA_summer} and \ref{fig:Power_SA_winter}). Even higher utilization factors are not economically optimal since a significantly higher storage capacity would be needed to completely follow the power peaks. The maximum power that the Power-to-Methanol plant can consume will therefore be lower than the installed capacity of the renewable park as already observed in other works \cite{Gorre2020}\cite{Chen2021}.

Storage is essential in the stand-alone configuration to grant continuous operation of the methanol synthesis plant in both scenarios, especially in low power production hours. 
\ce{H2} storage always plays a relevant role, though its optimal volume gets lower when allowing increased flexibility: the optimizer chooses larger methanol synthesis plant sizes rather than large storage capacities (Table \ref{tab:Optimization_SA}).
The battery instead plays a minor role compared to the \ce{H2} storage for the considered scenarios: the maximum amount of stored energy can satisfy the maximum power demand of the electrolyzer for less than 1 hour (if flexibility is allowed). Nevertheless, it intensively exchanges power  with the system especially during the production peaks and shortages when the electrolyzer cannot be operated (Figures \ref{fig:Power_SA_summer} and \ref{fig:Power_SA_winter}) and when the electrolyzer cannot be operated at full capacity (also see Figures A.\ref{fig:Storage_SA_summer} and A.\ref{fig:Storage_SA_winter}).


\begin{table}[h!bt]
    \centering
    \caption{Optimization values for the stand-alone case study for the considered scenarios (May-Jun and Nov-Dec). Refer to Appendix \ref{Additional_results} for the complete results.}
    \begin{tabular}{lccccccc}
        \hline
         &\textbf{Flexibility}& \makecell{\textbf{Specific cost} \\ \textbf{of MeOH}}& \textbf{$\mathrm{En_{batt}}$}& \textbf{$V$}& \textbf{$N_{\mathrm{mod}}$}& \textbf{$\dot{m}_{\mathrm{H_{2,MeOH},nom}}$} & \textbf{$\mathrm{MeOH_{y}}$} \\
         &\%/h & \texteuro/kg & MWh & \si{m^{3}} &   - & t/h & kt\\
                    \hline
        May-Jun& 0 &  1.78  & 65  & 5000* & 20 & 0.67  & 28.2 \\
        May-Jun& 5 & 1.40   & 27  & 800  & 33 & 1.43  & 37.9 \\
        May-Jun & 25  & 1.39 & 28  & 510  & 33 & 1.43  & 37.9 \\
                  \hline
       Nov-Dec & 0 & 2.45  & 42  & 5000*  & 14 & 0.44 & 18.4 \\
       Nov-Dec &5  &  1.70  & 23 & 940  & 28 & 1.22  & 29.2 \\
       Nov-Dec &25 & 1.69  & 17  & 900  & 26 & 1.15  & 28.5 \\
         \hline
         * bound \\
         \hline
    \end{tabular}    
    \label{tab:Optimization_SA}
\end{table}
\subsection{Discussion}
The optimal design of the Power-to-Methanol plant with and without storage in both case studies is highly dependent on the considered scenario. Nevertheless, general trends can be observed.

Flexible operation turned out to be an effective way to reduce the methanol production cost in Power-to-Methanol plants both with and without storage. The cost reduction due to the flexibility for both the grid-connected and stand-alone case studies is in line with other works in literature \cite{Gorre2020}\cite{Chen2021} even though a quantitative comparison is difficult as it depends on the considered scenarios and boundaries.

For the grid-connected case study, the reduction of the methanol production during the high-price hours allows reducing the operating costs and, at the same time, the production cost of methanol (value in agreement with the cost range calculated by IRENA (0.8-1.6 \$/kg) \cite{IRENA2021}). However, the optimal capacity factor of the methanol synthesis plant highly depends on the electricity price scenario (94 \% and 60 \% for the investigated scenarios).

In the stand-alone case study, flexible operation allows increasing the utilization of the electricity from the renewable park and the produced methanol, thus reducing the production cost. A similar cost reduction between flexible and not flexible stand-alone Power-to-Methanol plants can be found in literature (1.50-1.85 \$/kg and 1.89-2.84 \$/kg for flexible and not flexible Power-to-Methanol with hydrogen storage \cite{Chen2021}). Further improvements could be achieved in case some electricity is taken from another source \cite{Chen2021} or the plant is interfaced with the grid since electricity can be sold (instead of curtailed) or purchased in case of renewable generation surplus and shortage, respectively. 

The relative improvements in the objective function value get lower by increasing the maximum allowed ramp rate of the methanol synthesis plant. This aspect might be due to the considered time discretization since electricity price and renewable production may in practice vary on a shorter time scale.  
Nevertheless, flexibilization raises the issue of degradation: it might be worth not carrying out extremely quick load variations if these degrade the components of the plant, e.g., the catalyst, and the related economic savings are small.
The electrolyzers could particularly suffer from degradation due to frequent load changes, start-ups, and shut-downs. However, these phenomena are still not well understood and are difficult to model. An additional over-sizing of the electrolyzer to consider the long-term efficiency loss \cite{Sollai2023} or considering different operating strategies, e.g., avoiding frequent start-ups and shut-downs by introducing a penalty function, could be alternative ways to take into account these phenomena and refine the results. 

The role of storage depends on the considered boundaries, especially when the plant purchases electricity from the grid. When the electricity is cheap and its profile is not highly fluctuating as in the historical scenario 2021, the savings that the storage could allow are not sufficient to repay the investment. If the electricity price is higher as in the historical scenario 2022, storage improves the economic profitability of the plant. Nevertheless, flexibilization of the methanol synthesis plant has a greater economic benefit compared to adding storage.

For the considered boundary conditions, the \ce{H2} vessel provides big economic savings. The potential of the \ce{H2} vessel estimated in this optimization problem could even be rather conservative. In fact, a hydrogen bypass could be applied when the hydrogen consumption of the methanol synthesis plant is higher than the produced hydrogen ($\dot{m}_{\mathrm{H_{2},\, MeOH}}(t)\ge \dot{m}_{\mathrm{H_{2},\, prod}}(t)$): the produced hydrogen is compressed and directly fed into the methanol synthesis plant together with the already pressurized hydrogen taken from the vessel. The bypass would save electricity in the compression unit since the produced hydrogen has to be compressed up to the pressure of the synthesis plant (75 bar) instead of the higher pressure inside the vessel. However, modeling this configuration would have introduced additional variables and complexity to the optimization problem and it was out of the scope of this work. 

Compared to the vessel, the benefit of battery on the production cost is lower for the considered boundary conditions despite the nearly ideal model; in fact, small capacities are chosen by the optimizer. Nevertheless, it plays an essential role in stand-alone configurations to supply electricity in case of shortage. 
Also, the battery is used for short-term storage (a few hours). For longer storage periods, the \ce{H2} vessel seems to be a better option since it does not suffer from self-discharge. Nevertheless, the potential of batteries could change with other boundary conditions. In fact, their cost is expected to significantly decrease in the upcoming years \cite{IEA2021}. Furthermore, batteries could help smoothen the load changes of the plant and have other revenues derived from ancillary services to the grid, e.g., primary reserve and black start assistance, which have not been considered in this work. 

Finally, the importance of the considered scenarios and boundaries has to be once more highlighted. As shown for the investigated case studies and scenarios, the optimal design could significantly change with respect to the considered inputs. If there is no certainty about the scenarios, multiple scenarios could be considered in a two-stage stochastic design and scheduling optimization problem. 
Also, the design of the plant for the stand-alone case study could be optimized over longer time frames to better account for the seasonal behavior of renewable generation. However, reformulations of the optimization problem or time aggregation series techniques \cite{Schafer2020} should be used to reduce the computational burden.
\section{Conclusions} \label{Conclusions}
A combined design and scheduling optimization problem was formulated in GAMS for a Power-to-Methanol plant with both a battery and a hydrogen vessel to investigate the role of flexibility and storage. This problem was solved with the optimizer BARON for a grid-connected and stand-alone case study for different electricity price and renewable generation profiles.

Although overall global optimality cannot be guaranteed, the results suggest that flexible operation of the entire Power-to-Methanol plant significantly reduces the production cost of methanol in every case study. Also, it was noted that no significant cost reduction was achieved by considering maximum ramp rates higher than 10 \%/h for the considered scenarios and time discretization.

In the grid-connected case study, storage technologies play a relevant role when the electricity profile is highly fluctuating and at a high average price. In this scenario, the Power-to-Methanol configuration with both the battery and the \ce{H2} storage minimizes the methanol production cost. However, among these two storage technologies, the \ce{H2} storage provides the major economic benefits since it allows the decoupling between hydrogen production and conversion and thus a downsizing of the methanol synthesis plant.
In the stand-alone case study, storage is essential in all the scenarios, and it is widely used regardless of whether the methanol synthesis plant is operated constantly at nominal load or flexibly. 

Furthermore, it must be mentioned that longer time series should be considered for the optimal plant design to account for renewable production seasonality since the optimal designs differ when considering different generation scenarios. However, some time-aggregation series techniques \cite{Schafer2020}, multi-period optimization \cite{Martin2016}, or problem reformulations should be applied to numerically handle these optimization problems \cite{Mitsos2018}.
Finally, some models of the units could be further improved, e.g., by considering operating maps or degradation phenomena to increase the accuracy of the results. Also, the built optimization framework can be adapted to test other configurations, e.g., with the hydrogen bypass, or used for other Power-to-X processes by adjusting the key parameters.
\newpage
\textbf{Acknowledgments}\\
The authors gratefully acknowledge the financial support by the German Federal Ministry of E-ducation and Research (BMBF) within the H2Giga project DERIEL (grant number 03HY122D). The authors would also like to thank Dominik P. Goldstein for the creation of a submodel inserted in the Aspen Plus flowsheet of the methanol synthesis plant.\\
\textbf{CRediT authorship contribution statement}\\
\textbf{Simone Mucci}: Conceptualization, Methodology, Software, Investigation, Visualization, Writing - Original Draft. \textbf{Alexander Mitsos}: Conceptualization, Methodology, Supervision, Writing - Review \& Editing, Funding acquisition. \textbf{Dominik Bongartz}: Conceptualization, Methodology, Supervision, Writing - Review \& Editing, Funding acquisition, Project administration.

    \nocite{*}
    \printbibliography[heading=bibintoc] 

@book{IRENA2021,
address = {Abu Dhabi},
author = {{International Renewable Energy Agency (IRENA) and Methanol Institute}},
isbn = {9789292603205},
%pages = {1-124},
title = {{Innovation outlook: Renewable methanol}},
year = {2021}
}

@article{Ulonska2018,
title = {{Optimization of multiproduct biorefinery processes under consideration of biomass supply chain management and market developments}},
author = {Ulonska, Kirsten and K{\"{o}}nig, Andrea and Klatt, Marten and Mitsos, Alexander and Viell, J{\"{o}}rn},
doi = {10.1021/acs.iecr.8b00245},
journal = {Industrial and Engineering Chemistry Research},
number = {20},
pages = {6980--6991},
volume = {57},
year = {2018}
}

@article{Khadraoui2022b,
title = {{Production of lignin-containing cellulose nanofibrils by the combination of different mechanical processes}},
author = {Khadraoui, Malek and Khiari, Ramzi and Bergaoui, Latifa and Mauret, Evelyne},
doi = {10.1016/j.indcrop.2022.114991},
journal = {Industrial Crops and Products},
pages = {114991},
volume = {183},
year = {2022}
}

@article{Assen2016,
title = {{Selecting \ce{CO2} sources for \ce{CO2} utilization by environmental-merit-order curves}},
author = {{Von Der Assen}, Niklas and M{\"{u}}ller, Leonard J. and Steingrube, Annette and Voll, Philip and Bardow, Andr{\'{e}}},
doi = {10.1021/acs.est.5b03474},
issn = {15205851},
journal = {Environmental Science and Technology},
number = {3},
pages = {1093--1101},
volume = {50},
year = {2016}
}

@techreport{InternationalEnergyAgencyIEA2022,
address = {Paris},
author = {{International Energy Agency (IEA)}},
%pages = {1--159},
title = {{Renewables 2022: Analysis and forecast to 2027}},
url = {https://www.iea.org/reports/renewables-2022},
year = {2022}
}

@article{Burre2020,
title = {{Power-to-X: Between electricity storage, e-production, and demand side management}},
author = {Burre, Jannik and Bongartz, Dominik and Br{\'{e}}e, Luisa and Roh, Kosan and Mitsos, Alexander},
doi = {10.1002/cite.201900102},
issn = {15222640},
journal = {Chemie-Ingenieur-Technik},
number = {1-2},
pages = {74--84},
publisher = {Wiley-VCH Verlag},
volume = {92},
year = {2020}
}

@article{Zhang2016,
   author = {Qi Zhang and Ignacio E. Grossmann},
   doi = {10.1016/j.cherd.2016.10.006},
   issn = {02638762},
   journal = {Chemical Engineering Research and Design},
   month = {12},
   pages = {114-131},
   publisher = {Institution of Chemical Engineers},
   title = {Enterprise-wide optimization for industrial demand side management: Fundamentals, advances, and perspectives},
   volume = {116},
   year = {2016},
}

@article{Roh2019,
author = {Roh, Kosan and Br{\'{e}}e, Luisa C. and Perrey, Karen and Bulan, Andreas and Mitsos, Alexander},
journal = {Computer Aided Chemical Engineering},
doi = {10.1016/B978-0-12-818634-3.50296-4},
issn = {15707946},
pages = {1771--1776},
title = {Optimal Oversizing and Operation of the Switchable Chlor-Alkali Electrolyzer for Demand Side Management},
volume = {46},
year = {2019}
}

@article{Wang2020,
   author = {Ganzhou Wang and Alexander Mitsos and Wolfgang Marquardt},
   doi = {10.1002/aic.16947},
   issn = {15475905},
   issue = {6},
   journal = {AIChE Journal},
   pages = {1-9},
   title = {Renewable production of ammonia and nitric acid},
   volume = {66},
   year = {2020}
}

@article{Valverde2016,
author = {Valverde, Luis and Bordons, Carlos and Rosa, Felipe},
doi = {10.1109/TIE.2015.2465355},
issn = {02780046},
journal = {IEEE Transactions on Industrial Electronics},
number = {1},
pages = {167--177},
publisher = {Institute of Electrical and Electronics Engineers Inc.},
title = {{Integration of fuel cell technologies in renewable-energy-based microgrids optimizing operational costs and durability}},
volume = {63},
year = {2016}
}

@article{Allman2019,
author = {Allman, Andrew and Palys, Matthew J. and Daoutidis, Prodromos},
doi = {10.1002/aic.16434},
issn = {15475905},
journal = {AIChE Journal},
number = {7},
title = {{Scheduling-informed optimal design of systems with time-varying operation: A wind-powered ammonia case study}},
volume = {65},
year = {2019}
}

@article{SchulteBeerbuhl2015,
author = {{Schulte Beerb{\"{u}}hl}, S. and Fr{\"{o}}hling, M. and Schultmann, F.},
doi = {10.1016/j.ejor.2014.08.039},
issn = {03772217},
journal = {European Journal of Operational Research},
number = {3},
pages = {851--862},
title = {{Combined scheduling and capacity planning of electricity-based ammonia production to integrate renewable energies}},
volume = {241},
year = {2015}
}

@article{Allman2018,
author = {Allman, Andrew and Daoutidis, Prodromos},
doi = {10.1016/j.cherd.2017.10.010},
issn = {02638762},
journal = {Chemical Engineering Research and Design},
pages = {5--15},
publisher = {Institution of Chemical Engineers},
title = {{Optimal scheduling for wind-powered ammonia generation: Effects of key design parameters}},
volume = {131},
year = {2018}
}

@article{Chen2021,
   author = {Chao Chen and Aidong Yang},
   doi = {10.1016/j.enconman.2020.113673},
   OPTissn = {01968904},
   journal = {Energy Conversion and Management},
   OPTmonth = {1},
   pages = {113673},
   title = {Power-to-Methanol: The role of process flexibility in the integration of variable renewable energy into chemical production},
   volume = {228},
   year = {2021},
}

@article{Osman2020,
   author = {Ola Osman and Sgouris Sgouridis and Andrei Sleptchenko},
   doi = {10.1016/j.jclepro.2020.121627},
   OPTissn = {09596526},
   journal = {Journal of Cleaner Production},
   pages = {121627},
   title = {Scaling the production of renewable ammonia: A techno-economic optimization applied in regions with high insolation},
   volume = {271},
   year = {2020},
}

@inproceedings{Li2020a,
   author = {Jiarong Li and Jin Lin and Yonghua Song},
   doi = {10.1109/PESGM41954.2020.9282084},
   OPTisbn = {978-1-7281-5508-1},
   OPTissn = {19449933},
   booktitle = {2020 IEEE Power \& Energy Society General Meeting (PESGM)},
   pages = {1-5},
   publisher = {IEEE},
   title = {Capacity optimization of hydrogen buffer tanks in renewable Power to Ammonia ({P2A}) system},
   year = {2020},
}

@article{Gorre2019,
title = {Production costs for synthetic methane in 2030 and 2050 of an optimized Power-to-Gas plant with intermediate hydrogen storage},
journal = {Applied Energy},
volume = {253},
pages = {113594},
year = {2019},
issn = {0306-2619},
doi = {https://doi.org/10.1016/j.apenergy.2019.113594},
author = {Jachin Gorre and Felix Ortloff and Charlotte {van Leeuwen}},
}

@article{Gorre2020,
   author = {Jachin Gorre and Fabian Ruoss and Hannu Karjunen and Johannes Schaffert and Tero Tynjälä},
   doi = {10.1016/j.apenergy.2019.113967},
   OPTissn = {03062619},
   OPTissue = {October},
   journal = {Applied Energy},
   pages = {113967},
   title = {Cost benefits of optimizing hydrogen storage and methanation capacities for Power-to-Gas plants in dynamic operation},
   volume = {257},
   year = {2020},
}

@online{GAMS, 
    title = {{General Algebraic Modeling System (GAMS)}},
    OPTyear = 2023,
    url = {https://www.gams.com/},
    urldate = {2023-05-19}
}

@online{BARON, 
%    title = {{BARON Solver}},
%    OPTyear = 2023,
%    url = {https://www.minlp.com/baron-solver},
%    urldate = {2023-03-29}
%}

@article{BARON_2018,
author ={Aida Khajavirad and Nikolaos V. Sahinidis},
title={A hybrid LP/NLP paradigm for global optimization relaxations},
journal = {Mathematical Programming Computation},
volume={10},
pages={383-421},
year = {2018},
doi = {https://doi.org/10.1007/s12532-018-0138-5}
}

@article{Mitsos2018,
author = {Mitsos, Alexander and Asprion, Norbert and Floudas, Christodoulos A. and Bortz, Michael and Baldea, Michael and Bonvin, Dominique and Caspari, Adrian and Sch{\"{a}}fer, Pascal},
doi = {10.1016/j.compchemeng.2018.03.013},
issn = {00981354},
journal = {Computers \& Chemical Engineering},
pages = {209--221},
title = {{Challenges in process optimization for new feedstocks and energy sources}},
volume = {113},
year = {2018}
}

@article{Mucci2023,
title = {Power-to-X processes based on PEM water electrolyzers: A review of process integration and flexible operation},
journal = {Computers \& Chemical Engineering},
volume = {175},
pages = {108260},
year = {2023},
doi = {https://doi.org/10.1016/j.compchemeng.2023.108260},
author = {Simone Mucci and Alexander Mitsos and Dominik Bongartz}
}

@article{Schafer2020,
author = {Sch{\"{a}}fer, Pascal and Schweidtmann, Artur M. and Lenz, Philipp H.A. and Markgraf, Hannah M.C. and Mitsos, Alexander},
doi = {10.1016/j.compchemeng.2019.106598},
issn = {00981354},
journal = {Computers \& Chemical Engineering},
title = {{Wavelet-based grid-adaptation for nonlinear scheduling subject to time-variable electricity prices}},
volume = {132},
year = {2020}
}

@book{Tremel2018,
author = {Tremel, Alexander},
doi = {10.1007/978-3-319-72459-1},
isbn = {978-3-319-72458-4},
%pages = {1--100},
publisher = {Springer Cham},
series = {SpringerBriefs in Applied Sciences and Technology},
title = {{Electricity-based Fuels}},
year = {2018}
}

@book{IEA2021,
title = {{World energy outlook 2021}},
author = {{International Energy Agency (IEA)}},
pages = {1--386},
year = {2021}
}

@book{IRENA2017,
   author = {{International Renewable Energy Agency (IRENA)}},
   isbn = {978-92-9260-038-9},
   issue = {October},
   pages = {1-132},
   title = {Electricity storage and renewables: Costs and markets to 2030},
   year = {2017}
}

@article{Gonzalez-Castellanos2020,
title = {Detailed Li-ion battery characterization model for economic operation},
journal = {International Journal of Electrical Power \& Energy Systems},
volume = {116},
pages = {105561},
year = {2020},
issn = {0142-0615},
doi = {https://doi.org/10.1016/j.ijepes.2019.105561},
author = {Alvaro Gonzalez-Castellanos and David Pozo and Aldo Bischi},
}

@article{Martin2016,
author = {Mart{\'{i}}n, Mariano},
doi = {10.1016/j.compchemeng.2016.05.001},
issn = {00981354},
journal = {Computers \& Chemical Engineering},
pages = {43--54},
publisher = {Elsevier Ltd},
title = {{Methodology for solar and wind energy chemical storage facilities design under uncertainty: Methanol production from \ce{CO2} and hydrogen}},
volume = {92},
year = {2016}
}

@article{Carmo2013,
title = {{A comprehensive review on PEM water electrolysis}},
author = {Carmo, Marcelo and Fritz, David L. and Mergel, J{\"{u}}rgen and Stolten, Detlef},
journal = {International Journal of Hydrogen Energy},
OPTdoi = {10.1016/j.ijhydene.2013.01.151},
OPTkeywords = {Electrocatalysts, Electrolyzer separator plates, Hydrogen economy,PEM electrolysis modeling,PEM electrolyzers,Proton exchange membrane},
number = {12},
pages = {4901--4934},
volume = {38},
year = {2013}
}

@article{Jaervinen2022,
title = {{Automized parametrization of PEM and alkaline water electrolyzer polarisation curves}},
author = {J{\"{a}}rvinen, Lauri and Puranen, Pietari and Kosonen, Antti and Ruuskanen, Vesa and Ahola, Jero and Kauranen, Pertti and Hehemann, Michael},
OPTdoi = {10.1016/j.ijhydene.2022.07.085},
OPTissn = {03603199},
journal = {International Journal of Hydrogen Energy},
number = {75},
pages = {31985--32003},
volume = {47},
year = {2022}
}

@techreport{IEA2019,
address = {Paris},
author = {{International Energy Agency (IEA)}},
%pages = {1--203},
title = {{The future of hydrogen: Seizing today´s opportunities}},
url = {https://doi.org/10.1787/1e0514c4-en},
year = {2019}
}

@article{Kopp2017,
   author = {M. Kopp and D. Coleman and C. Stiller and K. Scheffer and J. Aichinger and B. Scheppat},
   OPTdoi = {10.1016/j.ijhydene.2016.12.145},
   issn = {03603199},
   issue = {19},
   journal = {International Journal of Hydrogen Energy},
   pages = {13311-13320},
   title = {Energiepark {Mainz}: Technical and economic analysis of the worldwide largest Power-to-Gas plant with PEM electrolysis},
   volume = {42},
   year = {2017},
}

@article{Tahan2022,
title = {Recent advances in hydrogen compressors for use in large-scale renewable energy integration},
journal = {International Journal of Hydrogen Energy},
volume = {47},
number = {83},
pages = {35275-35292},
year = {2022},
issn = {0360-3199},
doi = {https://doi.org/10.1016/j.ijhydene.2022.08.128},
author = {Mohammad-Reza Tahan},
}

@article{He2018,
title = {Compression performance optimization considering variable charge pressure in an adiabatic compressed air energy storage system},
journal = {Energy},
volume = {165},
pages = {349-359},
year = {2018},
issn = {0360-5442},
doi = {https://doi.org/10.1016/j.energy.2018.09.168},
author = {Yang He and Haisheng Chen and Yujie Xu and Jianqiang Deng},
}

@article{Frank2022,
title = {Loss Characterization of Advanced {VIGV} Configurations With Adjustable Blade Geometry},
journal = {Journal of Turbomachinery},
volume = {144},
number = {3},
pages = {031012},
year = {2022},
doi = {https://doi.org/10.1115/1.4052409},
author = {Frank, Roman G. and Wacker, Christian  and  Niehuis Reinhard},
}

@article{Moradi2019,
title = {{Hydrogen storage and delivery: Review of the state of the art technologies and risk and reliability analysis}},
author = {Moradi, Ramin and Groth, Katrina M.},
doi = {10.1016/j.ijhydene.2019.03.041},
issn = {03603199},
journal = {International Journal of Hydrogen Energy},
number = {23},
pages = {12254--12269},
publisher = {Elsevier Ltd},
url = {https://doi.org/10.1016/j.ijhydene.2019.03.041},
volume = {44},
year = {2019}
}

@article{Elberry2021,
title = {{Large-scale compressed hydrogen storage as part of renewable electricity storage systems}},
author = {Elberry, Ahmed M. and Thakur, Jagruti and Santasalo-Aarnio, Annukka and Larmi, Martti},
doi = {10.1016/j.ijhydene.2021.02.080},
issn = {03603199},
journal = {International Journal of Hydrogen Energy},
number = {29},
pages = {15671--15690},
volume = {46},
year = {2021}
}

@online{H2_storage_cost, 
    author = {{U.S. Department of energy (DOE)}},
    title = {Technical Targets for Hydrogen Delivery},
    OPTyear = 2022,
    url = {https://www.energy.gov/eere/fuelcells/doe-technical-targets-hydrogen-delivery},
    urldate = {2023-05-19}
}

@book{Hirscher2010,
address = {Weinheim, Germany},
author = {Hirscher, Michael},
doi = {10.1002/9783527629800},
isbn = {9783527629800},
pages = {1--365},
publisher = {Wiley-VCH Verlag GmbH \& Co. KGaA},
title = {{Handbook of Hydrogen Storage}},
year = {2010}
}

@article{Sollai2023,
   author = {Stefano Sollai and Andrea Porcu and Vittorio Tola and Francesca Ferrara and Alberto Pettinau},
   doi = {10.1016/j.jcou.2022.102345},
   issn = {22129820},
   journal = {Journal of \ce{CO2} Utilization},
   OPTmonth = {2},
   title = {Renewable methanol production from green hydrogen and captured \ce{CO2}: A techno-economic assessment},
   volume = {68},
   year = {2023},
}

@article{Anicic2014,
title = {{Comparison between two methods of methanol production from carbon dioxide}},
author = {Anicic, B., and Trop, P. and Goricanec, D.},
journal = {Energy},
OPTdoi = {10.1016/j.energy.2014.09.069},
pages = {279--289},
volume = {77},
year = {2014}
}

@article{Van-Dal2013,
title = {{Design and simulation of a methanol production plant from \ce{CO2} hydrogenation}},
author = {Van-Dal, {\'{E}}verton Sim{\~{o}}es and Bouallou, Chakib},
OPTdoi = {10.1016/j.jclepro.2013.06.008},
journal = {Journal of Cleaner Production},
OPTkeywords = {Aspen,CO2 mitigation,Methanol synthesis,Plus,Synthetic fuel},
pages = {38--45},
volume = {57},
year = {2013}
}

@article{VandenBussche1996,
title = {{A steady-state kinetic model for methanol synthesis and the water gas shift reaction on a commercial \ce{Cu/ZnO/Al_{2}O_{3}}}},
author = {Vanden Bussche, K. M. and Froment, G. F.},
journal = {Journal of Catalysis},
number = {161},
pages = {1--10},
year = {1996}
}

@article{Dieterich2020,
   author = {Vincent Dieterich and Alexander Buttler and Andreas Hanel and Hartmut Spliethoff and Sebastian Fendt},
   OPTdoi = {10.1039/d0ee01187h},
   issue = {10},
   journal = {Energy and Environmental Science},
   pages = {3207-3252},
   title = {Power-to-liquid via synthesis of methanol, DME or Fischer–Tropsch-fuels: a review},
   volume = {13},
   year = {2020},
}

@online{FlexMeOH22,
    author = {{BSE engineering}},
    title = {Power-to-Methanol at Small-Scale: FlexMethanol 10 {MW} \& 20 {MW} Module},
    url = {http://www.wfbe.de/BSE-Flyer-BASF-DB_web.pdf},
    OPTyear = 2023,
    urldate = {2023-05-19}
}

@article{Ramirez-Santos2017,
   author = {Álvaro A. Ramírez-Santos and Christophe Castel and Eric Favre},
   doi = {10.1016/j.memsci.2016.12.033},
   issn = {18733123},
   journal = {Journal of Membrane Science},
   pages = {191-204},
   publisher = {Elsevier},
   title = {Utilization of blast furnace flue gas: Opportunities and challenges for polymeric membrane gas separation processes},
   volume = {526},
   year = {2017}
}

@online{ENTSOE,
    author = {{ENTSO-E transparency platform}},
    title = {ENTSO-E: Day ahead prices},
    url = {https://transparency.entsoe.eu/transmission-domain/r2/dayAheadPrices/show?name=&defaultValue=false&viewType=TABLE&areaType=BZN&atch=false&dateTime.dateTime=23.11.2022+00:00|CET|DAY&biddingZone.values=CTY|10Y1001A1001A83F!BZN|10Y1001A1001A82H&resolution.values=PT60M&dateTime.timezone=CET_CEST&dateTime.timezone_input=CET+(UTC+1)+/+CEST+(UTC+2)},
    OPTyear = 2022,
    urldate = {2022-11-23}
}

@article{Reksten2022,
title = {{Projecting the future cost of PEM and alkaline water electrolysers; a CAPEX model including electrolyser plant size and technology development}},
author = {Reksten, Anita H. and Thomassen, Magnus S. and M{\o}ller-Holst, Steffen and Sundseth, Kyrre},
doi = {10.1016/j.ijhydene.2022.08.306},
issn = {03603199},
journal = {International Journal of Hydrogen Energy},
year = {2022},
issue = {90},
pages = {38106--38113},
volume = {47}
}

@book{Biegler1997,
title = {Systematic Methods of Chemical Process Design},
chapter={4-5},
author = {Lorenz T. Biegler and Ignacio E. Grossmann and Arthur W. Westerberg},
ISBN= {9780134924229},
Editor= {{Prentice Hall PTR}},
pages = {110--174},
year = {1997}
}

@article{Guthrie1969,
title = {{Data and techniques for preliminary capital cost estimating}},
author = {K.M. Guthrie},
issue = {76},
journal = {Chemical Engineering},
pages = {114-142},
volume = {3},
year = {1969}
}

@online{Methanolplants_cost,
   author = {{ADI Analytics LLC, Houston, Texas}},
   pages = {1-22},
   title = {Natural gas utilization via small scale methanol technologies},
   url = {http://www.sgicc.org/uploads/8/4/3/1/8431164/bftp_methanol_white_paper_vf.pdf},
   year = {2015},
   urldate = {2023-05-19}
}

@book{IRENA_RenGen_2021,
   author = {{International Renewable Energy Agency (IRENA)}},
   isbn = {978-92-9260-244-4},
   issn = {1476-4687},
   city = {Abu Dhabi},
   pages = {1-204},
   pmid = {25246403},
   title = {Renewable power generation costs in 2021},
   url = {https://www.irena.org/-/media/Files/IRENA/Agency/Publication/2018/Jan/IRENA_2017_Power_Costs_2018.pdf},
   year = {2022},
}
\newpage
\appendix
\section{Appendix} \label{Appendix}
\captionsetup{labelformat=AppendixTables} 
\setcounter{table}{0} 
\setcounter{figure}{0} 

\subsection{PEM water electrolyzer efficiency} \label{PEMWE_efficiency}
The implemented efficiency law of the PEM water electrolyzer ($\mathrm{\eta_{PEM}}$) is a function of the power of the electrolysis module ($P_{\mathrm{mod}}$) and the operating pressure $p_{\mathrm{PEM}}$, and it is defined as follows:

\begin{equation*} \mathrm{\eta_{PEM}}=a_{00}+a_{10}\cdot P_{\mathrm{mod}}+a_{20}\cdot P_{\mathrm{mod}}^2+a_{01}\cdot p_{\mathrm{PEM}}, \end{equation*}

where $P_{\mathrm{mod}}$ is equal to the overall power of the electrolyzer $P_{\mathrm{PEM}}$ divided by the number of modules ($N_{\mathrm{mod}}$).

The values of the parameters were obtained by fitting the results of the electrolyzer model \cite{Jaervinen2022} in the power and pressure ranges of 0.2-2.5 MW and 20-40 bar, respectively, and are collected in Table A.\ref{tab:Efficiency_PEMWE_parameters}.

\begin{table}[h!bt]
    \centering
    \caption{Parameters for the efficiency law of the electrolyzer ($\mathrm{R^2}$=0.999).}
    \begin{tabular}{lcc}
        \hline
         & \textbf{Values} & \textbf{Unit}\\
          \hline
        $a_{00}$ &  $0.813$ & - \\
        $a_{10}$ &  $-1.010 \cdot 10^{-01}$ & \si{MW^{-1}}\\
        $a_{20}$ &  $+1.397 \cdot 10^{-02}$ & \si{MW^{-2}}\\
        $a_{01}$ &  $-3.118 \cdot 10^{-04} $ & \si{bar^{-1}} \\
         \hline
    \end{tabular}    
    \label{tab:Efficiency_PEMWE_parameters}
\end{table}
\subsection{Thermophysical properties}  \label{Properties}
The thermophysical properties of hydrogen from the Aspen database (Peng-Robinson property method) were used to estimate the parameters of $c_{\mathrm{p}} $ and $k$, and as the benchmark for the model results.

Table A\ref{tab:Compr} shows the comparison between the calculated power of the compressor unit and the Aspen Plus data for a hydrogen inlet stream of 1900 kg/h at 40 bar for different pressure ratios ($\beta$). The error for the considered operating points is lower than 1 \%. 

\begin{table}[h!bt]
    \centering
    \caption{Compressor model validation.}
    \begin{tabular}{ccc}
        \hline
        $\beta$ & \textbf{Model} & \textbf{Aspen data}\\
          \hline
        2.0 &  710 kW  & 710 kW \\
        2.5 &  970 kW & 973 kW\\
        3.0 &  1196 kW & 1202 kW\\
        3.5 &  1397 kW & 1406 kW\\
         \hline
    \end{tabular}    
    \label{tab:Compr}
\end{table}
In Figure A\ref{fig:H2_density}, the density of hydrogen at high pressure and its linear approximation for the hydrogen storage model is shown.
\begin{figure}[h!bt]
    \centering
    \includegraphics[width=0.5\textwidth]{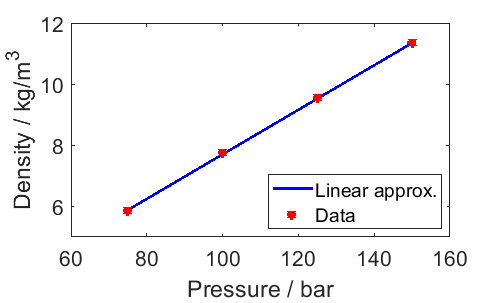}    
    \caption{\ce{H2} density at 25 °C as function of the pressure (data from the Aspen database).}
    \label{fig:H2_density}
\end{figure}
\subsection{Hydrogen storage cost} \label{CAPEX_vessel}
The considered hydrogen storage cost is 500 \$/kg, which corresponds to the target cost for 160-bar stationary gaseous hydrogen tanks for the year 2020 \cite{H2_storage_cost}.
This value was converted from the specific cost for unit of hydrogen mass into the specific cost for unit of volume and updated from year 2007 to year 2021 as follows:

\begin{equation*} \mathrm{Cost_{vessel}= 500 \; \frac{\$}{kg}\cdot 0.073\; \frac{kg}{m^3 bar} \cdot (160\; bar -1\; bar) \cdot UF} , \end{equation*}

where $0.073\; \mathrm{\frac{kg}{m^3 bar} \cdot (160\; bar -1\; bar)}$ determines the available hydrogen mass content per unit of volume, and $\mathrm{UF}$ is the update factor equal to 1.35. The calculated cost is around 7800  \$/\si{m^{3}}. Lower specific costs for vessels might be possible in scenarios with higher production volumes \cite{H2_storage_cost} and considering a lower maximum pressure inside the vessel, which reduces the wall thickness, thus the material cost.
\subsection{Methanol synthesis plant}   \label{Methanol}
The plant was modeled in Aspen Plus \textsuperscript{\textregistered} V11 (Figure A.\ref{fig:MeOH_flowsheet}). The property methods \textit{PSRK} and \textit{NRTL} were used for the thermodynamic properties in high and low pressure units, respectively.
\begin{figure}[h!bt]
    \centering
    \includegraphics[angle=90, height=0.96\textheight]{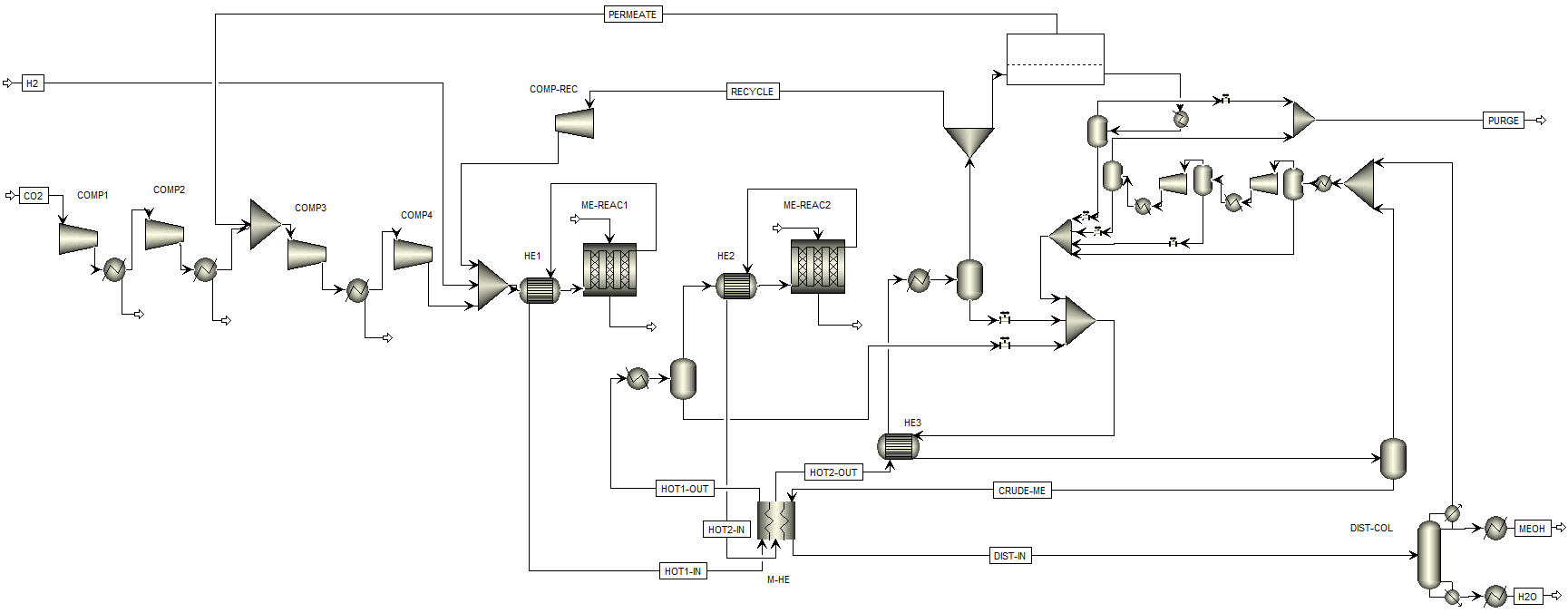}    
    \caption{Methanol synthesis plant model in Aspen Plus\textsuperscript{\textregistered} V11.}
    \label{fig:MeOH_flowsheet}
\end{figure}
The methanol synthesis occurs in a two-stage multi-tubular reactor, which is cooled via evaporating water at around 245 °C. The operating pressure of the reactor unit is around 75 bar, and a 1 bar pressure drop was assumed for each reactor stage. A commercial catalyst Cu/ZnO/\ce{Al2O3} is considered for the implemented kinetic model LHHW (Langmuir Hinshelwood Hougen Watson) \cite{Van-Dal2013}\cite{VandenBussche1996}.
The choice of having two reactor stages is motivated by the low conversion per pass. The intermediate removal of products via condensation helps shift the chemical equilibrium. 

After the second stage of the reactor, the unreacted gases are separated from crude methanol, the mixture of water and methanol, via condensation and partially (98 \%) recycled back to the reactor. Also, a membrane separation process \cite{Ramirez-Santos2017} was considered to partially (90 \%) recover \ce{H2} from the purge stream.

Crude methanol is then purified via distillation to meet the desired purity specification, i.e., 99.85 \%wt. (AA grade methanol). Heat for the reboiler can be supplied via heat integration. In particular, the surplus heat of the methanol reactor and the heat from the purge stream combustion are sufficient to satisfy its thermal demand. However, this heat integration is not shown in the flowsheet in Figure A.\ref{fig:MeOH_flowsheet}, though it was considered for the calculation of the heat and cooling demands. 

Table A.\ref{tab:MeOH_ref_data} summarizes the key information about the methanol synthesis plant that was used as the reference for the technical and economical analysis.

\begin{table}[h!bt]
    \centering
    \caption{Main data of the reference methanol synthesis plant.}
    \begin{tabular}{lc}
        \hline
         & \textbf{Values} \\
          \hline
        \ce{H2} mass flow rate &  1.9 t/h \\
        \ce{CO2} mass flow rate &  15.4 t/h \\
        \ce{CH3OH} mass flow rate &  9.9 t/h \\
        Overall electricity demand &  1.8 MW\\
        Overall cooling demand &  13.8 MW\\
        Net heat demand &  0 MW\\
        First law efficiency &  84.3 \% \\
         \hline
    \end{tabular}    
    \label{tab:MeOH_ref_data}
\end{table}

\FloatBarrier
\subsection{Methanol synthesis plant cost} \label{CAPEX_MeOH}
The CAPEX of the main equipment units of the methanol synthesis plant (Fig. A. \ref{fig:MeOH_flowsheet}), i.e., compressors, reactor, heat exchangers, flashes, and distillation column, was estimated by using the cost models proposed by Biegler et al. \cite{Biegler1997}. The cost of the membrane was estimated according to Ramírez-Santos et al.'s model \cite{Ramirez-Santos2017}. 

Biegler et al.'s method \cite{Biegler1997} is based on the well-established Guthrie's method \cite{Guthrie1969}, which estimates the bare module cost of the equipment with accuracy within the range of $\pm$ 25-40 \%. A multiplying factor of 1.85 \cite{Biegler1997} was then used to account for fixed and working capital costs. An update factor from 1968 to 2021 of 6.15 was considered. 

The estimated CAPEX of the reference methanol synthesis plant with a yearly production of 0.08 Mt is 27.8 M\texteuro ~(350 \texteuro/(t/y)). This value is within the typical range for large-scale methanol synthesis plants (200-700 \$/(t/y) \cite{Methanolplants_cost}), which means that the capital cost is probably underestimated. On a similar production capacity (0.06 Mt/y), the estimated investment cost of a small-scale methanol synthesis plant in Texas was around 40 M\$ \cite{Methanolplants_cost}. The difference might be due to the uncertainty of the cost models and the different boundary conditions (the steam methane reforming and hydrogen compression units are also included in the plant in Texas). The estimated investment cost is in line with an e-methanol plant in Norway (yearly production of 0.1 Mt and a capital cost of 200 M\$ \cite{IRENA2021}) if the cost of the electrolysis unit and its replacement after 10 years are considered.
\subsection{Scenarios} \label{Scenarios_Appendix}
In this section, additional details about the electricity price and renewable availability scenarios are shown. Table A.\ref{tab:Variations_table} shows that a relevant part of the hourly variations of the electricity price and renewable power availability profiles (relative to the mean value over the considered scenario) are lower than 5\%. This explains why moderate flexibility of the methanol synthesis plant (in \%/h) is already highly beneficial from an economic perspective. However, a linear correlation between the variability of the scenarios and the flexibility of the methanol synthesis plant cannot be established since other factors, e.g., the storage utilization, operating constraints, and chosen objective function play a role.

\begin{table}[h!bt]
    \centering
    \caption{Percentage of the hourly price or renewable power availability variations between consecutive hours lower than 5\%, 10\%, and 25\% for the considered scenarios. The relative hourly variations were calculated with respect to the average value over the considered scenario, i.e., 51 \texteuro/MWh, 154 \texteuro/MWh, 52 MW, and 39 MW for the `Historical 2021', `Historical 2022', `May-Jun', and `Nov-Dec' scenarios, respectively.}
    \begin{tabular}{ccccc}
        \hline
         \textbf{Hourly variations} & \textbf{Historical 2021} & \textbf{Historical 2022}& \textbf{May-Jun}& \textbf{Nov-Dec}\\
          \hline
         $\leq$ 5\% &  48 \%  & 42 \%  &36 \%  &52 \%  \\
         $\leq$ 10\% &  72 \%  & 63 \%  &50 \%  &75 \%  \\
         $\leq$ 25\% &  96 \%  & 93 \%  &80 \%  &97 \%  \\
         \hline
    \end{tabular}    
    \label{tab:Variations_table}
\end{table}
\subsection{Optimization results} \label{Additional_results}
In this section, additional optimization results are provided.
\subsubsection{Grid-connected case study}
The optimal design variables and the key performance metrics for the two scenarios of the grid-connected case study are shown in Tables A.\ref{tab:Optimization_GC_2021_complete} and A.\ref{tab:Optimization_GC_2022_complete}.

Figure A.\ref{fig:Storage_utilization_GC_Scen2} shows how the battery and the \ce{H2} storage are optimally used for the high-price electricity scenario. As expected, the energy stored in the battery and the hydrogen in the vessel are consumed in a medium-to-short term, thus confirming that the chosen time frame for scheduling (2 months) is suitable.
While the use of the battery is similar between the `No flexibility' and `25 \%/h flexibility' cases, major differences can be noticed in the \ce{H2} storage. The amount of hydrogen stored inside the vessel is significantly lower when the methanol synthesis plant is operated flexibly. A similar qualitative behavior can instead be noticed in the pressure profile (dependent on both the amount of hydrogen and the storage volume). The differences could be explained by the different optimal operating strategies.

\begin{figure}[h!bt]
    \centering
    \includegraphics[width=1\textwidth]{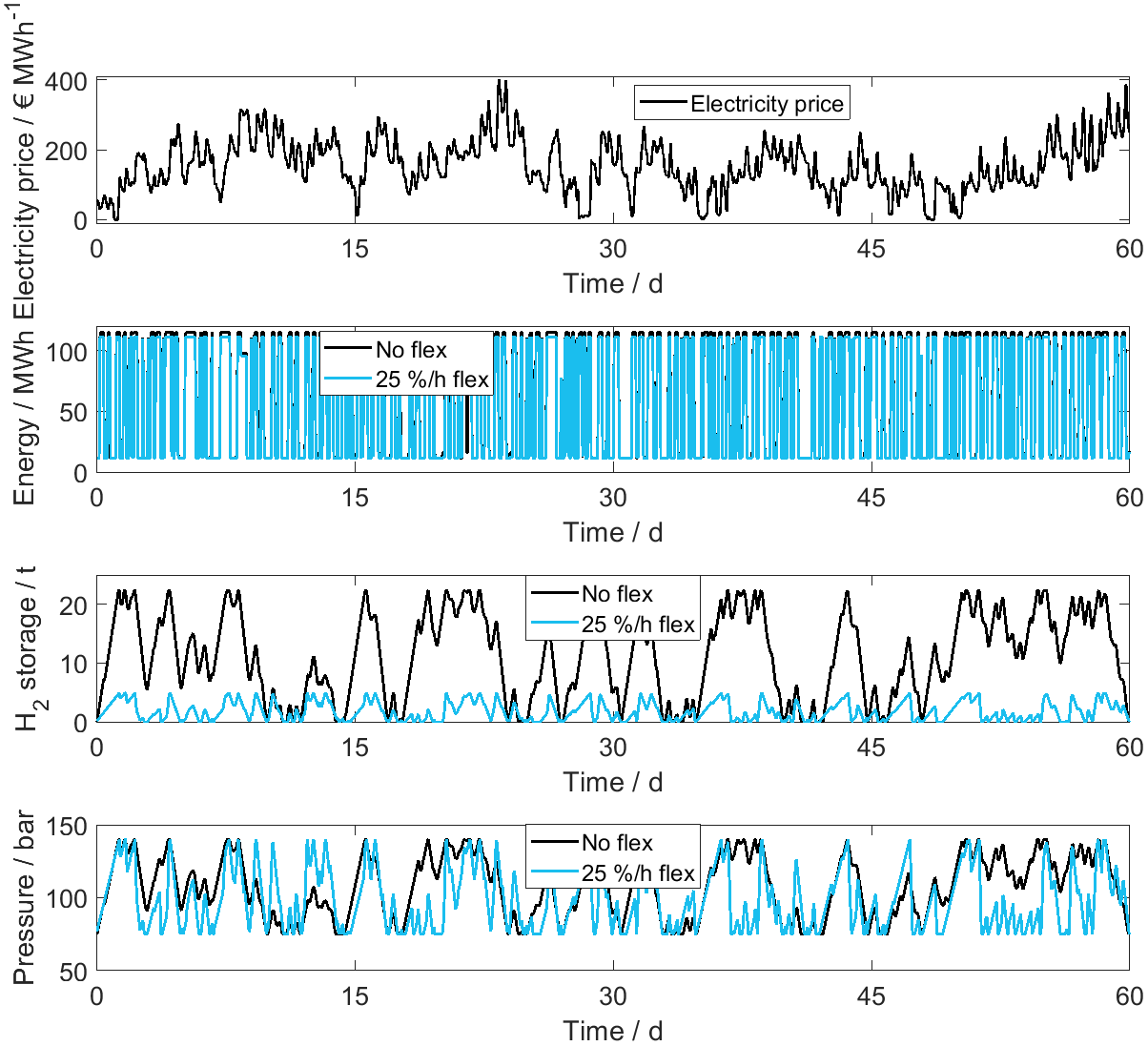}  
    \caption{Energy content of the battery, and \ce{H2} content and pressure inside the \ce{H2} vessel for the `No flexibility' and `25 \%/h flexibility' cases. Only the results of the `Battery \& Vessel' configuration are shown. The electricity cost profile is also shown to better interpret the optimal scheduling. Note: the energy profiles of the battery for the `No flexibility' and `25 \%/h flexibility' cases are similar, thus they can hardly be distinguished in the picture. }
    \label{fig:Storage_utilization_GC_Scen2}
\end{figure}

\begin{table}[h!bt]
    \centering
    \caption{Optimization results for the grid-connected case study for the low-price low-variable electricity price scenario for the 4 plant configurations, i.e., with both storage (B\&V), with the hydrogen vessel (V), with battery (B), and without any storage (No B\&V).}
\rotatebox{90}{%
    \begin{tabular}{lccccccccccc}
        \hline
        &\textbf{Flexibility}& \makecell{\textbf{Specific cost} \\ \textbf{of MeOH}} & \textbf{$\mathrm{En_{batt}}$} & \textbf{$V$} & \textbf{$N_{\mathrm{mod}}$} & \textbf{$p_{\mathrm{PEM}}$} &  \textbf{$\beta_{\mathrm{max}}$} & \textbf{$\dot{m}_{\mathrm{H_{2,MeOH},nom}}$} & \textbf{$\mathrm{MeOH_{y}}$} & \textbf{$\mathrm{CAPEX}_{0}$} & \textbf{$\mathrm{OPEX}_{y}$}\\
        & \%/h&\texteuro/kg & MWh & \si{m^{3}} & - & bar & - &  t/h &   kt &  M\texteuro & M\texteuro\\
        \hline
        B\&V& 0 & 0.91  & 5*  & 25* & 40* & 40*  & 1.88* & 1.91  &  80.3   & 135  & 51 \\
        V & 0& 0.91 & - & 25* & 40* & 40*  & 1.88*& 1.91  & 80.3   & 133  & 51 \\
        B &0&  0.91   & 5*  & - & 40* & 40*  & 1.88*& 1.91& 80.3   & 134  & 51 \\
        No B\&V & 0 & 0.91   & - & - & 40* & 40*  & 1.88*& 1.91  &  80.3   & 133  & 51 \\
        \hline
        B\&V& 5& 0.90  & 5*  & 25*  & 40* & 40*  & 1.88*& 1.91  &  76.3  & 135  & 46 \\
        V &5& 0.90  & - & 25* & 40* & 40*  & 1.88*& 1.91  & 76.2   & 133  & 46\\
        B &5 & 0.91  & 5* & - & 40* & 40*  & 1.88* & 1.91 & 76.3   & 134  & 46 \\
        No B\&V&5 & 0.90  & - & - & 40* & 40*  & 1.88* & 1.91  & 76.2   & 133  & 46 \\
        \hline
        B\&V& 10& 0.90  & 5*  & 25*  & 40* & 40*  & 1.88* & 1.91  &  75.7  & 135  & 46 \\
        V &10& 0.90  & - & 25*  & 40* & 40*  & 1.88* & 1.91  &  75.7  & 133  & 46\\
        B&10 & 0.90  & 5*  & - & 40* & 40* & 1.88* & 1.91 & 75.7   & 134 & 46 \\
        No B\&V&10 & 0.90  & - & - & 40* & 40*& 1.88* & 1.91  & 75.7  & 133  & 46 \\
        \hline
        * bound\\
        \hline
       \end{tabular} 
}%
    \label{tab:Optimization_GC_2021_complete}
\end{table}

\begin{table}[h!bt]
    \centering
    \caption{Optimization results for the grid-connected case study for the high-price high-variable electricity price scenario for the 4 plant configurations, i.e., with both storage (B\&V), with the hydrogen vessel (V), with battery (B), and without any storage (No B\&V).}
    \rotatebox{90}{%
    \begin{tabular}{lccccccccccc}
        \hline
        &\textbf{Flexibility}& \makecell{\textbf{Specific cost} \\ \textbf{of MeOH}} & \textbf{$\mathrm{En_{batt}}$} & \textbf{$V$} & \textbf{$N_{\mathrm{mod}}$} & \textbf{$p_{\mathrm{PEM}}$} &  \textbf{$\beta_{\mathrm{max}}$} & \textbf{$\dot{m}_{\mathrm{H_{2,MeOH},nom}}$} & \textbf{$\mathrm{MeOH_{y}}$} & \textbf{$\mathrm{CAPEX}_{0}$} & \textbf{$\mathrm{OPEX}_{y}$}\\
                & \%/h&\texteuro/kg & MWh & \si{m^{3}} & - & bar & - &  t/h &   kt &  M\texteuro & M\texteuro\\
          \hline
        B\&V& 0&  1.86  & 115  & 4730  & 40* & 40*  & 3.5* & 1.11 &  47.0   & 196  & 57 \\
        V &0& 1.88  & - & 4940  & 40* & 40* & 3.5* & 1.15 &  48.5   & 162 & 65 \\
        B & 0& 2.01  & 116 & - & 40* & 31.6  & 2.38 & 1.91  &  80.7   & 171  & 135 \\
        No B\&V &0& 2.02 & - & - & 40* & 31.5 & 2.38 & 1.91  &  80.7  & 135  & 141 \\
        \hline
        B\&V& 5& 1.69  & 116  & 1210 & 40* & 40*  & 3.5* & 1.74 &  42.9  & 180  & 44 \\
        V&5 & 1.70  & - & 1320  & 40* & 40*  & 3.5* & 1.74  &  42.7   & 144  & 49 \\
        B&5 & 1.74 & 115 & - & 40* & 40* & 1.88* & 1.91  &  49.2   & 169 & 59 \\
        No B\&V &5& 1.75  & - & - & 40* & 40*  & 1.88* & 1.91  &  49.5  & 133 & 64 \\
        \hline
        B\&V& 10&1.68  & 116 & 1130 & 40* & 40* & 3.5* & 1.74 &  42.4   & 179  & 43 \\
        V & 10&1.69  & - & 1250 & 40* & 40*  & 3.5* & 1.74 &  42.7  & 144 & 49 \\
        B&10 & 1.72 & 111 & - & 40* & 40*  & 1.88* & 1.91  &  49.3  & 168  & 58 \\
       No B\&V &10& 1.73   & - & - & 40* & 40*  & 1.88* & 1.91 &  49.5  & 133  & 64 \\
        \hline
       B\&V&25& 1.68  & 111  & 1040  & 40* & 40* & 3.5* & 1.75  &  42.9   & 177  & 44 \\
        V &25& 1.69 & - & 1100  & 40* & 40*  & 3.5* & 1.74&   42.7 & 143 & 49 \\
        B&25 & 1.71  & 115  & - & 40* & 40*  & 1.88* & 1.91  &  49.3  & 169  & 57 \\
        No B\&V &25& 1.72 & - & - & 40* & 40*  & 1.88* & 1.91  &  49.4   & 133  & 63 \\
        \hline
             * bound \\
             \hline
    \end{tabular}  
    }%
    \label{tab:Optimization_GC_2022_complete}
\end{table}
\FloatBarrier
\subsubsection{Stand-alone case study}
The optimal design variables and the key performance metrics for the two scenarios of the stand-alone case study are shown in Table A.\ref{tab:Optimization_SA_2022_complete}.

Figures A.\ref{fig:Storage_SA_summer} and A.\ref{fig:Storage_SA_winter} show how the battery and the \ce{H2} storage are optimally operated for the considered scenarios. 
Both storage technologies are operated on a daily basis to cope with the high power peaks due to the PV generation during the May-Jun scenario (Figure A.\ref{fig:Storage_SA_summer}) when flexibility of the methanol synthesis plant is allowed. In fact, the possibility of varying the methanol production rate in periods of renewable production shortage reduces the need for high amounts of stored hydrogen and long-term storage.
A similar operation can be noticed for the Nov-Dec scenario (Figure A.\ref{fig:Storage_SA_winter}). The main difference lies in the longer duration of the `hold phase' of both storage technologies, and it is probably due to the scenario, which is characterized by a relevant wind power generation that has different characteristic times from PV generation. Nevertheless, these hold phases last a few days at most since the methanol synthesis plant is assumed to be always operating. 
Thus, on the one hand, shorter time frames could have been considered for scheduling; on the other hand, even longer time frames would allow for taking into account the seasonality of power generation in the design phase. This aspect is particularly evident from the design differences when considering two power generation profiles of the same renewable park in two periods of the year (Table A.\ref{tab:Optimization_SA_2022_complete}). In fact, the optimal design of the Power-to-Methanol plant for the Nov-Dec scenario leads to a higher production cost of methanol when the plant is optimally operated in the May-Jun scenario (1.43 \texteuro/kg instead of 1.39 \texteuro/kg).

\begin{figure}[h!bt]
    \centering
    \includegraphics[width=1\textwidth]{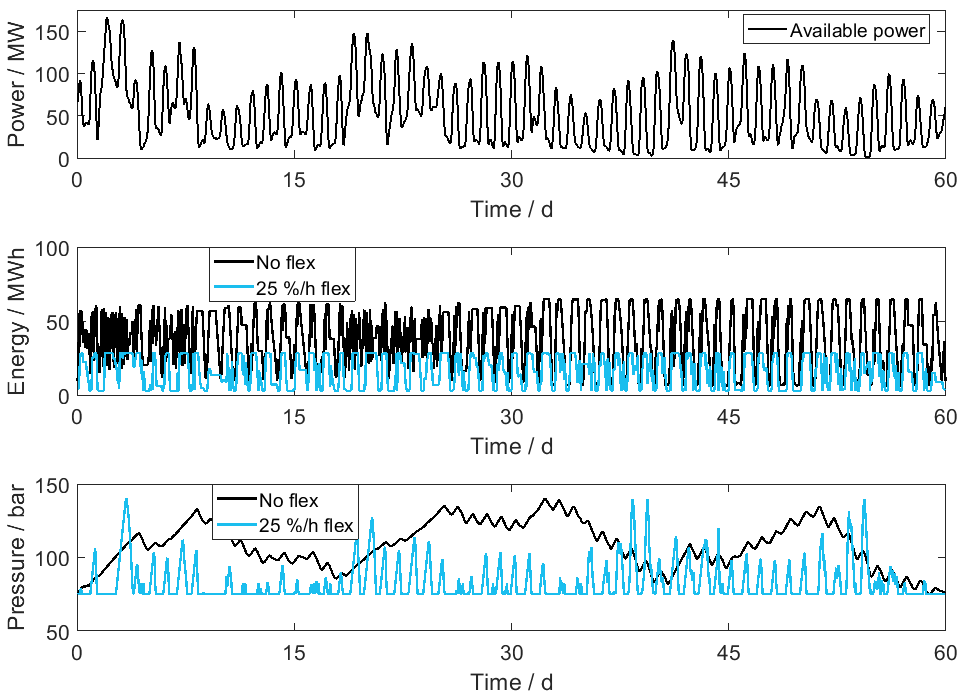}  
    \caption{Energy content of the battery and \ce{H2} pressure inside the \ce{H2} vessel for the `No flexibility' and `25 \%/h flexibility' cases. The power availability profile (May-Jun scenario) is also shown to better interpret the optimal scheduling.}
    \label{fig:Storage_SA_summer}
\end{figure}
\begin{figure}[h!bt]
    \centering
    \includegraphics[width=1\textwidth]{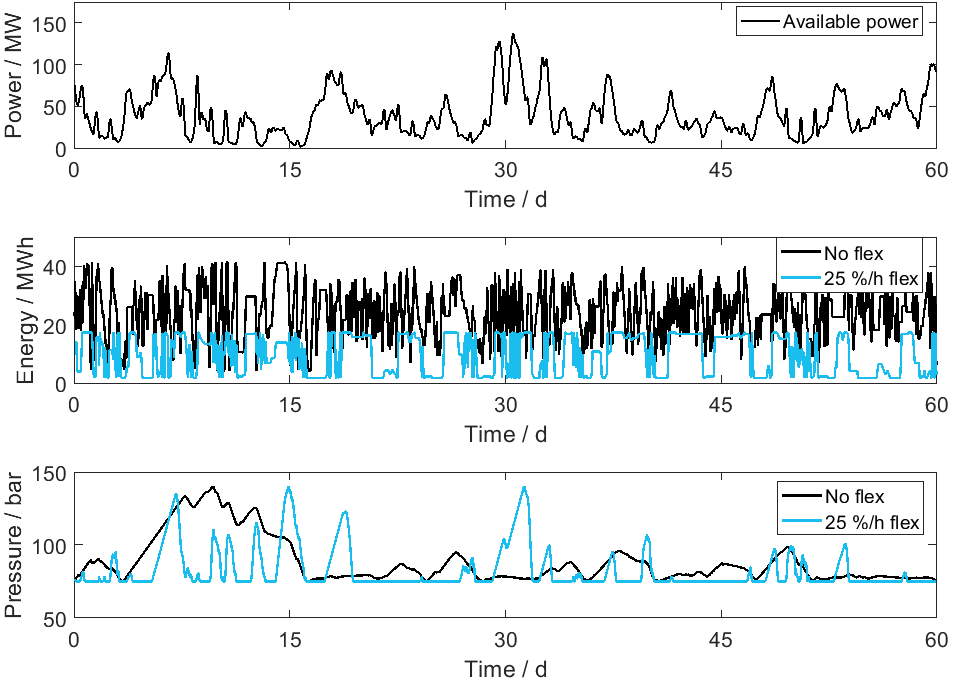}  
    \caption{Energy content of the battery and \ce{H2} pressure inside the \ce{H2} vessel for the `No flexibility' and `25 \%/h flexibility' cases. The power availability profile (Nov-Dec scenario) is also shown to better interpret the optimal scheduling.}
    \label{fig:Storage_SA_winter}
\end{figure}

\begin{table}[h!bt]
    \centering
    \caption{Optimization results for the stand-alone case study for the considered scenarios (May-Jun and Nov-Dec).}
    \rotatebox{90}{%
    \begin{tabular}{lcccccccccccc}
        \hline
        &\textbf{Flexibility}& \makecell{\textbf{Specific cost} \\ \textbf{of MeOH}} & \textbf{$\mathrm{En_{batt}}$} & \textbf{$V$} & \textbf{$N_{\mathrm{mod}}$} & \textbf{$p_{\mathrm{PEM}}$} &  \textbf{$\beta_{\mathrm{max}}$} & \textbf{$\dot{m}_{\mathrm{H_{2,MeOH},nom}}$} & \textbf{$\mathrm{MeOH_{y}}$} & \textbf{$\mathrm{CAPEX}_{0}$} & \textbf{$\mathrm{OPEX}_{y}$} & \textbf{$\mathrm{P}_{renew}$}\\
        &\%/h &   \texteuro/kg & MWh & \si{m^{3}} & - &  bar & - &  t/h &   kt & M\texteuro &  M\texteuro & \%\\
          \hline
        May-Jun&0&  1.78  & 65 & 5000*  & 20 & 40*  & 3.5* & 0.67  & 28.2   & 348  & 2  & 71 \%\\
        May-Jun&5& 1.40  & 27  & 800  & 33 & 40* & 3.5* &  1.43  &  37.9  & 350  & 3   & 93 \%\\
        May-Jun&25& 1.39 & 28 &  510  & 33 & 40*  & 3.5* & 1.43  & 37.9   & 349  & 3   & 93 \%\\
        \hline
        Nov-Dec&0&   2.45  & 42  & 5000*  & 14 & 40* & 3.5* & 0.44 & 18.4   & 322  & 1   & 62 \%\\
        Nov-Dec&5& 1.70  & 23  & 940  & 28 & 40*  & 3.5* &  1.22  & 29.2  & 335 & 2   & 95 \%\\
        Nov-Dec&25& 1.69  & 17   & 840  & 26 & 40*  & 3.5* & 1.15 &  28.5  & 327  & 2   & 93 \%\\
        \hline
        * bound \\
        \hline
    \end{tabular}  
    }%
    \label{tab:Optimization_SA_2022_complete}
\end{table}
\end{document}